\documentclass[a4,psfig,english,11pt,epsf,portrait]{article}
\usepackage{amsmath,amsfonts,latexsym,amscd,amssymb,theorem}
\usepackage{epsfig}
\usepackage{amsmath}
\usepackage{psfrag}
\def\wt{\widetilde}

\def\R{\hbox{\bf R}}
\def\Z{\hbox{\bf Z}}

\def\o{\overline}

\def\N{\hbox{\bf N}}

\def\g{\gamma}
\def\b{\beta}

\def\a{\alpha}
\def\h{\hat}

\def\intws{\int_{\mathbb{R}^{n+1}}}

\def\l{{\lambda}}

\def\ti{\tilde}

\def\inm{-\!\!\!\!\!\!\!\!\;\int}

\def\<{\langle}
\def\>{\rangle}
\newcommand{\ba}{\begin{eqnarray}}
\newcommand{\ea}{\end{eqnarray}}

\newtheorem{thm}{Theorem}[section]

\newtheorem{theorem}[thm]{Theorem}
\newtheorem{definition}[thm]{Definition}
\newtheorem{lemma}[thm]{Lemma}
\newtheorem{proposition}[thm]{Proposition}
\newtheorem{corollary}[thm]{Corollary}
\newtheorem{rem}[thm]{Remark}

\numberwithin{equation}{section}

\renewcommand{\R}{{\mathbb R}}
\renewcommand{\Z}{{\mathbb Z}}
\renewcommand{\N}{{\mathbb N}}

\renewcommand{\g}{\gamma}

\begin{document}

\title{\bf On a parabolic logarithmic Sobolev inequality}
\author{
\normalsize\textsc{ H. Ibrahim $^{*,\ddagger,\dagger}$, R. Monneau
  \footnote{Universit\'{e} Paris-Est, CERMICS,
Ecole des Ponts, 6 et 8 avenue Blaise Pascal, Cit\'e
Descartes Champs-sur-Marne, 77455 Marne-la-Vall\'ee Cedex 2, France.
E-mails: ibrahim@cermics.enpc.fr, monneau@cermics.enpc.fr
\newline \indent $\,\,{}^\dagger$CEREMADE, Universit\'{e}
Paris-Dauphine, Place De Lattre de Tassigny, 75775 Paris Cedex 16,
France
\newline \indent $\,\,{}^\ddagger$LaMA-Liban,
Lebanese University, P.O. Box 826 Tripoli, Lebanon}}}
\vspace{20pt}

\maketitle


\centerline{\small{\bf{Abstract}}} \noindent{\small{In order to extend
    the blow-up criterion of solutions to the Euler equations, Kozono
    and Taniuchi \cite{KT} have proved a logarithmic Sobolev inequality
    by means of isotropic (elliptic) $BMO$ norm. In this paper, we show
    a parabolic version of the Kozono-Taniuchi inequality by means of anisotropic
    (parabolic) $BMO$ norm. More precisely we give an upper bound for
    the $L^{\infty}$ norm of a function in terms of its parabolic $BMO$
    norm, up to a logarithmic correction involving its norm in some
    Sobolev space. As an application, we also explain how to apply this
    inequality in order to establish a long-time existence result for a
    class of nonlinear parabolic problems.}}

\hfill\break
 \noindent{\small{\bf{AMS subject classifications:}}} {\small{42B35,
     54C35, 42B25, 39B05.}}\hfill\break
  \noindent{\small{\bf{Key words:}}} {\small{Logarithmic
      Sobolev inequalities, parabolic $BMO$ spaces, anisotropic Lizorkin-Triebel
      spaces, harmonic analysis.}}\hfill\break

\section{Introduction and main results}\label{sec1}
In \cite{KT}, Kozono and Taniuchi showed an $L^{\infty}$ estimate of a given
function by means of its $BMO$ norm (space of functions of bounded mean
oscillation) and the logarithm of its
norm in some Sobolev space. In fact, they proved
that for $f\in W^{s}_{p}(\R^{n})$, $1<p<\infty$, the following estimate
holds (with $\log^{+}x = \max(\log x, 0)$):
\begin{equation}\label{koz_tan_bis}
\|f\|_{L^{\infty}(\R^{n})}\leq C(1 +
\|f\|_{BMO(\R^{n})}(1+\log^{+}\|f\|_{W^{s}_{p}(\R^{n})})),\quad sp>n,
\end{equation}
for some constant $C=C(n,p,s)>0$. The main advantage of the above
estimate is that it was successfully applied (see \cite[Theorem
2]{KT}) to extend the blow-up criterion of solutions to the Euler
equations which was originally given by Beale, Kato and Majda in
\cite{BKaM}. Inequality (\ref{koz_tan_bis}), as well as some
variants of it, are shown (see \cite{KT, Ogawa, KOT03}) using
harmonic analysis on isotropic functional spaces of the
Lizorkin-Triebel and Besov type. However, as is well known, it is
important, say for parabolic partial differential equations to
consider spaces that are anisotropic.

Motivated by the study of the long-time existence of a certain class
of singular parabolic coupled systems (see \cite{IJM_PI, CRAS_IJM}),
we show in this paper an analogue of the Kozono-Taniuchi inequality
(\ref{koz_tan_bis}) but of the parabolic (anisotropic) type. Due to
the parabolic anisotropy, we consider functional spaces on
$\R^{n+1}=\R^{n}\times \R$ with the generic variable $z=(x,t)$,
where each coordinate $x_{i}$, $i=1\cdots n$ is given the weight
$1$, while the time coordinate $t$ is given the weight $2$. We now
state the main results of this paper. The first result concerns a
Kozono-Taniuchi parabolic type inequality on the entire space
$\R^{n+1}$. Introducing parabolic bounded mean oscillation $BMO_{p}$
spaces, and parabolic Sobolev spaces $W_{2}^{2m,m}$ (for the
definition of these spaces, see Definitions \ref{para_bmo} and
{\ref{para_Sob}}), we present our first theorem.
\begin{theorem}{\bf (Parabolic logarithmic Sobolev inequality on
    $\R^{n+1}$)}\label{theorem1}\\
Let $u\in W_{2}^{2m,m}(\R^{n+1})$, $m>\frac{n+2}{4}$. Then there exists
a constant $C=C(m,n)>0$ such that:
\begin{equation}\label{cara_eq2}
\|u\|_{L^{\infty}(\R^{n+1})}\leq C \left(1+ \|u\|_{BMO_{p}(\R^{n+1})} \left(1+
    \log^{+} \|u\|_{W_{2}^{2m,m}(\R^{n+1})}\right) \right).
\end{equation}
\end{theorem}
The proof of Theorem \ref{theorem1} will be given in Section
\ref{sec2}, and is based on an approach developed by Ogawa
\cite{Ogawa}. Let us mention that our proof in this paper is self-contained. The
second result of this paper concerns a Kozono-Taniuchi parabolic type
inequality on the bounded domain
$$\Omega_{T}=(0,1)^{n}\times (0,T)\subset \R^{n+1},\quad T>0.$$
More precisely, our next theorem reads:
\begin{theorem}{\bf (Parabolic logarithmic Sobolev inequality on a bounded
    domain)}\label{theorem2}\\
Let $u\in W_{2}^{2m,m}(\Omega_{T})$ with $m>\frac{n+2}{4}$. Then there exists
a constant $C=C(m,n,T)>0$ such that:
\begin{equation}\label{theo2_eq1}
\|u\|_{L^{\infty}(\Omega_{T})}\leq C \left(1+
  \|u\|_{\overline{BMO}_{p}(\Omega_{T})} \left(1+
    \log^{+} \|u\|_{W_{2}^{2m,m}(\Omega_{T})}\right) \right),
\end{equation}
where $\|\cdot\|_{\overline{BMO}_{p}(\Omega_{T})} =
\|\cdot\|_{BMO_{p}(\Omega_{T})} + \|\cdot\|_{L^{1}(\Omega_T)}$. 
\end{theorem}
The proof of Theorem \ref{theorem2} will be
given in Section \ref{sec3}.
\subsection{Brief review of the literature}
The brief review presented here only concerns logarithmic Sobolev
inequalities of the elliptic type. Up to our knowledge, logarithmic
Sobolev inequalities of the parabolic type have not been treated
elsewhere in the literature.

The original type of the logarithmic
Sobolev inequalities was found in Brezis-Gallouet \cite{BreGal} and
Brezis-Wainger \cite{BreWai} where the authors
investigated the relation between $L^{\infty}$, $W^{k}_{r}$ and
$W^{s}_{p}$ and proved that there holds the embedding:
\begin{equation}\label{B_G_W}
\|f\|_{L^{\infty}(\R^{n})}\leq C
(1+\log^{\frac{r-1}{r}}(1+\|f\|_{W^{s}_{p}(\R^{n})})),\quad sp>n
\end{equation}
provided $\|f\|_{W^{k}_{r}(\R^{n})}\leq 1$ for $kr=n$. The estimate
(\ref{B_G_W}) was applied to prove global existence of solutions to the
nonlinear Schr\"{o}dinger equation (see \cite{BreGal, HayWol}). Similar
embedding for $f\in (W^{s}_{p}(\R^{n}))^{n}$ with $\mbox{div}f=0$
was investigated by Beale-Kato-Majda in \cite{BKaM}. The authors showed
that:
\begin{equation}\label{BKM_est}
\|\nabla f\|_{L^{\infty}}\leq C
(1+\|\mbox{rot}f\|_{L^{\infty}}(1+\log^{+}\|f\|_{W^{s+1}_{p}})+\|\mbox{rot}f\|_{L^{2}}),
\quad sp>n,
\end{equation}
where they made use of this estimate in order to give a blow-up
criterion of solutions to the Euler equations (see \cite{BKaM}). In \cite{KT}, Kozono and
Taniuchi showed their inequality (\ref{koz_tan_bis}) in order to extend the blow-up
criterion of solutions to the Euler equations given in
\cite{BKaM} (see \cite[Theorem 2]{KT}). A
generalized version of (\ref{koz_tan_bis}) in Besov spaces was given in
Kozono-Ogawa-Taniuchi \cite{KOT03}. Finally, a sharp version
of the logarithmic Sobolev inequality of the Beale-Kato-Majda and the
Kozono-Taniuchi type in the Lizorkin-Triebel spaces was established by
Ogawa in \cite{Ogawa}.
\subsection{Organization of the paper} This paper is organized as
follows. In Section \ref{sec2}, we recall basic
tools used in our analysis, and give the proof of Theorem
\ref{theorem1}. In Section \ref{sec3}, we present the proof of Theorem
\ref{theorem2}, and as an application, we explain how to use the
parabolic Kozono-Taniuchi inequality in order to prove the long-time
existence of certain parabolic equations.

\section{A parabolic Kozono-Taniuchi inequality on $\R^{n+1}$}\label{sec2}
This section is devoted to the proof of Theorem \ref{theorem1}.
\subsection{Preliminaries and basic tools}
\subsubsection{Parabolic $BMO_p$ and Sobolev spaces}
We start by recalling some definitions and introducing some notations. A
generic point in $\R^{n+1}$ will be denoted by $z=(x,t)\in
\R^{n}\times \R$, $x=(x_{1},.\,.\,.\,,\,x_{n})$. Let $S(\R^{n+1})$ be the
usual Schwartz space, and
$S'(\R^{n+1})$ the corresponding dual space. Let $u\in S'(\R^{n+1})$.
For $\xi=(\xi_{1},.\,.\,.\,,\,\xi_{n})\in \R^{n}$ and $\tau\in \R$
we denote by $\mathcal{F}u(\xi,\tau)\equiv\h{u}(\xi,\tau)$, and
$\mathcal{F}^{-1}u(\xi,\tau)\equiv\check{u}(\xi,\tau)$
the Fourier, and the inverse Fourier
transform of $u$ respectively. We also denote
$D_{t}^{r}=\frac{\partial^{r}}{\partial t^{r}}$,
$r\in \N$, and $D_{x}^{s}$, $s\in \N$, any derivative
with respect to $x$ of order $s$. The parabolic distance from $z=(x,t)$
to the origin is defined by:
\begin{equation}\label{para_distance}
\|z\|=\max\left\{|x_{1}|,.\,.\,.\,,\,|x_{n}|, |t|^{1/2}\right\}.
\end{equation}
Let $\mathcal{O}\subseteq \R^{n+1}$ be an open set.
The parabolic bounded mean oscillation space $BMO_p$ and the parabolic
Sobolev space $W_{2}^{2m,m}$ are now recalled.
\begin{definition}{\bf (Parabolic bounded mean oscillation spaces)}\label{para_bmo}\\
A function $u\in L^{1}_{loc}(\mathcal{O})$ is said to be of parabolic
bounded mean oscillation, $u\in BMO_{p}(\mathcal{O})$ if we have:
\begin{equation}\label{bmo_norm}
\|u\|_{BMO_{p}(\mathcal{O})}=\sup_{Q\subset \mathcal{O}}\frac{1}{|Q|}\int_{Q}|u-u_{Q}|<+\infty.
\end{equation}
Here $Q$ denotes an arbitrary parabolic cube
\begin{equation}\label{haAnas}
Q=Q_{r}=Q_{r}(z_{0})=\{z\in \R^{n+1};\, \|z-z_{0}\|<r\},
\end{equation}
and
\begin{equation}\label{essjfa}
u_{Q}=\frac{1}{|Q|}\int_{Q}u.
\end{equation}
The functions in $BMO_{p}$ are
defined up to an additive constant. We also define the space
$\overline{BMO}_{p}$ as:
$$\overline{BMO}_{p}(\mathcal{O})=BMO_{p}(\mathcal{O}) \cap
L^{1}(\mathcal{O})\quad \mbox{with}\quad \|\cdot\|_{\overline{BMO}_{p}} =
\|\cdot\|_{BMO_{p}} + \|\cdot\|_{L^{1}}.$$
\end{definition}
\begin{definition}{\bf (Parabolic Sobolev spaces)}\label{para_Sob}\\
Let $m$ be a non-negative integer. We define the parabolic Sobolev space
$W_{2}^{2m,m}(\mathcal{O})$ as follows:
$$W_{2}^{2m,m}(\mathcal{O})=\{u\in L^{2}(\mathcal{O}); D_{t}^{r}D_{x}^{s}u\in
L^{2}(\mathcal{O}),\, \forall r,s\in \N \mbox{ such that } 2r+s\leq 2m\}.$$
The norm of $u\in W_{2}^{2m,m}(\mathcal{O})$ is defined by:
$\|u\|_{W_{2}^{2m,m}(\mathcal{O})}=\sum^{2m}_{j=0}\sum_{2r+s=j}
\|D^{r}_{t}D^{s}_{x}u\|_{L^{2}(\mathcal{O})}$.
\end{definition}
The next lemma concerns a Sobolev embedding of $W_{2}^{2m,m}$.
\begin{lemma}{\bf (Sobolev embedding,  \cite[Lemma 3.3]{LSU})}\label{Sobo_emb}\\
Let $m$ be a non-negative integer satisfying $m> \frac{n+2}{4}$. Then
there exists a positive constant $C$ depending on $m$ and $n$ such that
for any $u\in
W_{2}^{2m,m}(\mathcal{O})$, the function $u$ is continuous and bounded on
$\mathcal{O}$, and satisfies
\begin{equation}\label{Embed_eq}
\|u\|_{L^{\infty}(\mathcal{O})}\leq C \|u\|_{W_{2}^{2m,m}(\mathcal{O})}.
\end{equation}
\end{lemma}
\subsubsection{Parabolic Lizorkin-Triebel and Besov spaces}
Here we give the definition of Lizorkin-Triebel spaces. These spaces are
constructed out of the parabolic Littlewood-Paley decomposition that we
recall here. Let $\psi_{0}(z)\in C_{0}^{\infty}(\R^{n+1})$ be a function
such that
\begin{equation}\label{psi_lp}
\psi_{0}(z)=1\,\mbox{ if }\, \|z\|\leq 1\quad \mbox{and}\quad
\psi_{0}(z)=0\,\mbox{ if }\, \|z\|\geq 2.
\end{equation}
For such a function $\psi_0$, we may define a smooth, anisotropic dyadic
partition of unity $(\psi_{j})_{j\in \N}$ by letting
$$\psi_{j}(z)=\psi_{0}(2^{-ja}z)-\psi_{0}(2^{-(j-1)a}z)\quad
\mbox{if}\quad j\geq 1.$$
Here $a=(1,.\,.\,.\,, 1,2)\in\R^{n+1}$, and for $\eta\in \R$,
$b=(b_{1},.\,.\,.\,, b_{n},b_{n+1})\in \R^{n+1}$, the dilation
$\eta^{b}z$ is defined by $\eta^{b}z=(\eta^{b_{1}}z_{1},.\,.\,.\,,
\eta^{b_{n}}z_{n},\eta^{b_{n+1}}z_{n+1})$. It is clear that
$$\sum^{\infty}_{j=0}\psi_{j}(z)=1\quad \mbox{for}\quad z\in \R^{n+1},$$
and
$$\mbox{supp}\,\psi_{j}\subset \{z;\, 2^{j-1}\leq \|z\|\leq
2^{j+1}\},\quad j\geq 1.$$ Define $\phi_{j}$, $j\geq 0$ as the
inverse Fourier transform of $\psi_{j}$, i.e.
$\hat{\phi}_{j}=\psi_{j}$. It is worth noticing that
\begin{equation}\label{hedshsu1}
\phi_{j}(z)=2^{(n+2)(j-1)}\phi_{1}(2^{(j-1)a}z)  \quad \mbox{for}\quad
j\geq 1,
\end{equation}
and that for any $u\in S'(\R^{n+1})$,
$$u=(2\pi)^{-\frac{(n+1)}{2}}\sum^{\infty}_{j=0}\phi_{j}*u,
\mbox{ with convergence in } S'(\R^{n+1}).$$ We now give the
definition of the anisotropic Besov and Lizorkin-Triebel spaces.
\begin{definition}{\bf (Anisotropic Besov spaces)}\label{aniso_besov}\\
The anisotropic Besov space $B_{p,q}^{s}(\R^{n+1})=B^{s}_{p,q}$, $s\in
\R$, $1\leq p\leq \infty$ and $1\leq q\leq \infty$ is the space of functions
$u\in S{'}(\R^{n+1})$ with finite quasi-norms
\begin{equation}\label{ani_bes_norm}
\|u\|_{B_{p,q}^{s}}=\left(\sum^{\infty}_{j=0}2^{sqj}\|\phi_{j}*u\|^{q}_{L^{p}(\R^{n+1})}\right)^{1/q}
\end{equation}
and the natural modification for $q=\infty$, i.e.
\begin{equation}\label{ani_bes_infini}
\|u\|_{B_{p,\infty}^{s}}=\sup_{j\geq 0} 2^{sj} \|\phi_{j}*u\|_{L^{p}(\R^{n+1})}.
\end{equation}
\end{definition}
\begin{definition}{\bf (Anisotropic Lizorkin-Triebel spaces)}\label{tri_lizo_def}\\
The anisotropic Lizorkin-Triebel space $F_{p,q}^{s}(\R^{n+1})=F^{s}_{p,q}$, $s\in
\R$, $1\leq p< \infty$ and $1\leq q\leq \infty$ (or $1\leq q<\infty$ and
$p=\infty$) is the space of functions
$u\in S{'}(\R^{n+1})$ with finite quasi-norms
\begin{equation}\label{ani_tri_lizo_norm}
\|u\|_{F_{p,q}^{s}}=\left\|\left(\sum^{\infty}_{j=0}2^{sqj}|\phi_{j}*u|^{q}\right)^{1/q}\right\|_{L^{p}(\R^{n+1})}
\end{equation}
and the natural modification for $q=\infty$, i.e.
\begin{equation}\label{ani_tri_lizo_infini}
\|u\|_{F_{p,\infty}^{s}}= \|\sup_{j\geq 0} 2^{sj} |\phi_{j}*u|
\|_{L^{p}(\R^{n+1})}.
\end{equation}
\end{definition}
A very useful space throughout our analysis will be the truncated
anisotropic (parabolic) Lizorkin-Triebel space
$\widetilde{F}^{s}_{p,q}$ that we define here.
\begin{definition}{\bf (Truncated anisotropic Lizorkin-Triebel space)}\label{trunT_L}\\
The truncated anisotropic Lizorkin-Triebel space
$\widetilde{F}_{p,q}^{s}(\R^{n+1})=\widetilde{F}^{s}_{p,q}$, $s\in
\R$, $1\leq p< \infty$ and $1\leq q\leq \infty$ ($1\leq q<\infty$ if
$p=\infty$) is the space of functions $u\in S{'}(\R^{n+1})$ with finite quasi-norms
\begin{equation}\label{trunT_Lnorm}
\|u\|_{\widetilde{F}_{p,q}^{s}}=\left\|\left(\sum^{\infty}_{j=1}2^{sqj}|
\phi_{j}*u|^{q}\right)^{1/q}\right\|_{L^{p}(\R^{n+1})}
\end{equation}
and the natural modification for $q=\infty$, i.e.
\begin{equation}\label{trunT_Lnorm_infi}
\|u\|_{\widetilde{F}_{p,\infty}^{s}}= \|\sup_{j\geq 1} 2^{sj} |\phi_{j}*u|
\|_{L^{p}(\R^{n+1})}.
\end{equation}
\end{definition}
The basic difference between $F^{s}_{p,q}$ and
$\widetilde{F}^{s}_{p,q}$ is that in $\widetilde{F}^{s}_{p,q}$ we
omit the term $\phi_{0}*u$ and only take in consideration the terms
$\phi_{j}*u$, $j\geq 1$. 
Sobolev embeddings of parabolic
Lizorkin-Triebel and Besov spaces are shown by the next two lemmas.
\begin{lemma}{\bf (Embeddings of Besov spaces,
    \cite[Theorem 7]{JohSic})}\label{lem1}\\
Let $s,t\in \R$, $s>t$, and $1\leq p,r\leq\infty$ satisfy:
$s-\frac{n+2}{p}=t-\frac{n+2}{r}$.
Then for any $1\leq q\leq \infty$ we have the following continuous embedding
\begin{equation}\label{johsic2}
B^{s}_{p,q}(\R^{n+1})\hookrightarrow B^{t}_{r,q}(\R^{n+1}).
\end{equation}
\end{lemma}
\begin{lemma}{\bf (Sobolev embeddings, \cite[Proposition
    2]{Stockert82})}\label{lem2}\\
Take an integer $m\geq 1$. Then we have
\begin{equation}\label{stock2}
B^{2m}_{2,1} \hookrightarrow W_{2}^{2m,m} \hookrightarrow
B^{2m}_{2,\infty}.
\end{equation}
\end{lemma}
\subsection{Basic logarithmic Sobolev inequality}
In this subsection we show a basic logarithmic Sobolev inequality. In
particular, we show the following lemma.
\begin{lemma}{\bf (Basic logarithmic Sobolev
    inequality)}\label{BSI_lelemme}\\
Let $u\in W_{2}^{2m,m}(\R^{n+1})$ for some $m\in \N$, $m>\frac{n+2}{4}$. Then
there exists some constant $C=C(m,n)>0$ such that
\begin{equation}\label{BSob_ineq}
\|u\|_{\wt{F}^{0}_{\infty,1}}\leq C\left(1+ \|u\|_{\wt{F}^{0}_{\infty,2}}\left(1+
(\log^{+}\|u\|_{W_{2}^{2m,m}})^{1/2}\right)\right).
\end{equation}
\end{lemma}
{\bf Proof.} First, let us mention that the ideas of the proof
of this lemma are inspired from the proof of Ogawa \cite[Corollary
2.4]{Ogawa}. The proof is divided into three steps, and the constants in
the proof may vary from line to line.\\

\noindent {\bf Step 1.} (Estimate of $\|u\|_{\wt{F}^{0}_{\infty,1}}$).\\

\noindent Let $\g>0$, and $N\in \N$ be two arbitrary
variables. We compute:
\begin{eqnarray*}
\|u\|_{\wt{F}^{0}_{\infty,1}}&\leq& \left\|\sum_{1\leq j<N}|\phi_{j}*u|\right\|_{L^{\infty}} +
\left\|\sum_{j\geq N} 2^{-\g j}2^{\g
    j}|\phi_{j}*u|\right\|_{L^{\infty}}\\
&\leq&N^{1/2}\left\|\left(\sum_{1\leq j<N}|\phi_{j}*u|^{2}
  \right)^{1/2}\right\|_{L^{\infty}} + C_{\g}2^{-\g N} \left\|\left(\sum_{j\geq N}
  (2^{\g j}|\phi_{j}*u| )^{2}\right)^{1/2} \right\|_{L^{\infty}}\\
&\leq& C_{\g}\left(N^{1/2} \|u\|_{\wt{F}^{0}_{\infty,2}} + 2^{-\g N}
\|u\|_{F^{\g}_{\infty,2}}\right),
\end{eqnarray*}
where $C_{\g}>0$ is a positive constant.\\

\noindent{\bf Step 2.} (Optimization in $N$).\\

\noindent We optimize the previous inequality in $N$ by setting:
$$N=1\quad \mbox{if}\quad \|u\|_{F^{\g}_{\infty,2}} \leq 2^{\g}
\|u\|_{\wt{F}^{0}_{\infty,2}}.$$
In this case we can easily check that:
\begin{equation}\label{Step1ineq}
\|u\|_{\wt{F}^{0}_{\infty,1}}\leq C_{\g}
\|u\|_{\wt{F}^{0}_{\infty,2}}\left(1+\left(\log^{+}\frac{\|u\|_{F^{\g}_{\infty,2}}}
{\|u\|_{\wt{F}^{0}_{\infty,2}}}\right)^{1/2}\right).
\end{equation}
In the case where $\|u\|_{F^{\g}_{\infty,2}} > 2^{\g}
\|u\|_{\wt{F}^{0}_{\infty,2}}$, we choose $1\leq \b < 2^{\g}$ such that
$$N=\log_{2^{\g}}^{+}\left(\b\frac{
    \|u\|_{F^{\g}_{\infty,2}}}{\|u\|_{\wt{F}^{0}_{\infty,2}}}\right)\in \N.$$
We then compute:
\begin{eqnarray*}
N^{1/2} \|u\|_{\wt{F}^{0}_{\infty,2}} + 2^{-\g N}
\|u\|_{F^{\g}_{\infty,2}}&\leq& \|u\|_{\wt{F}^{0}_{\infty,2}}
\left(\frac{1}{\b}+\left[\log^{+}_{2^{\g}}
\left(\b \frac{\|u\|_{F^{\g}_{\infty,2}}}{\|u\|_{\wt{F}^{0}_{\infty,2}}}\right)
\right]^{1/2} \right)\\
&\leq& \|u\|_{\wt{F}^{0}_{\infty,2}}
\left(\frac{1}{\b}+\left[\frac{2}{\log
    2^{\g}}\log^{+}\frac{\|u\|_{F^{\g}_{\infty,2}}}{
        \|u\|_{\wt{F}^{0}_{\infty,2}}} \right]^{1/2} \right)\\
&\leq& C_{\g} \|u\|_{\wt{F}^{0}_{\infty,2}} \left(1+ \left( \log^{+}\frac{\|u\|_{F^{\g}_{\infty,2}}}{
        \|u\|_{\wt{F}^{0}_{\infty,2}}}\right)^{1/2} \right),
\end{eqnarray*}
hence we also have (\ref{Step1ineq}) with a different constant
$C_{\g}$.\\

\noindent {\bf Step 3.} (Estimate of $\|u\|_{F^{\g}_{\infty,2}}$ and conclusion).\\

\noindent Noting the inequality
$$
x\left(\log \left(e+\frac{y}{x}\right) \right)^{1/2}\leq
\left\{
\begin{aligned}
& C\left(1+x (\log \left(e+y\right))^{1/2}\right)\quad &\mbox{for}& \quad 0<x\leq 1\\
& Cx(\log (e+y))^{1/2}\quad &\mbox{for}& \quad  x>1,
\end{aligned}
\right.
$$
we deduce from (\ref{Step1ineq}) that:
\begin{equation}\label{Step2ineq}
\|u\|_{\wt{F}^{0}_{\infty,1}}\leq C \left(1+\|u\|_{\wt{F}^{0}_{\infty,2}} \left(1+
  \left(\log^{+}\|u\|_{F^{\g}_{\infty,2}}  \right)^{1/2}\right)
\right),
\end{equation}
where the constant $C$ depends also on $\g$. We now estimate the term
$\|u\|_{F^{\g}_{\infty,2}}$. Choose $\g$ such that
$$0<\g< 2m-\frac{n+2}{2}.$$
Call $\a=2m-\frac{n+2}{2}$, we compute:
\begin{eqnarray}\label{slam}
\|u\|_{F^{\g}_{\infty,2}}&=& \left\|\left(\sum_{j\geq 0} 2^{2j\g}
    |\phi_{j}*u|^{2} \right)^{1/2} \right\|_{L^{\infty}}\nonumber\\
&\leq& \left(\sum_{j\geq 0}2^{2j (\g-\a)}\right)^{1/2}
\left\|\sup_{j\geq 0}2^{\a j} |\phi_{j}*u|\right\|_{L^{\infty}}\nonumber\\
&\leq & C \|u\|_{B^{\a}_{\infty,\infty}}.
\end{eqnarray}
It is easy to check (see (\ref{johsic2}), Lemma \ref{lem1}, and
(\ref{stock2}), Lemma \ref{lem2}) that we have the continuous embeddings
$$W_{2}^{2m,m}\hookrightarrow B^{2m}_{2,\infty}\hookrightarrow B^{\a}_{\infty,\infty}.$$
Therefore (from inequality (\ref{slam})) we get:
$$\|u\|_{F^{\g}_{\infty,2}}\leq C \|u\|_{W_{2}^{2m,m}},$$
hence the result directly follows from
(\ref{Step2ineq}). $\hfill{\blacksquare}$
\subsection{Proof of Theorem \ref{theorem1}}
In this subsection we present the proof of several lemmas leading to the
proof of Theorem \ref{theorem1}. We start with the following lemma concerning
mean estimates of functions on
parabolic cubes. Call $Q_{2^{j}}\subset \R^{n+1}$, $j\geq 0$, any arbitrary
parabolic cube of radius
$2^{j}$ (see (\ref{haAnas}) for the definition of parabolic cubes). For
the sake of simplicity, we denote
\begin{equation}\label{key_para}
Q^{j}=Q_{2^{j}}\quad \mbox{for all}\quad j\in \Z.
\end{equation}
Our next lemma reads:
\begin{lemma}{\bf (Mean estimates on parabolic cubes)}\label{mean_est}\\
Let $u\in BMO_{p}(\R^{n+1})$. Take
$Q^{j}\subset Q^{j+1}$, $j\geq 0$ ($Q^{j}$ and $Q^{j+1}$ do not
necessarily have the same center). Then we have (with the notation
(\ref{essjfa})): 
\begin{equation}\label{meanest_eq1}
|u_{Q^{j+1}}-u_{Q^{j}}|\leq
\left(1+2^{n+2}\right)\|u\|_{BMO_{p}}.
\end{equation}
More generally, we have for any $Q^{j}\subseteq Q^{k}$, $j,k\in \Z$:
\begin{equation}\label{EfS_eq1}
|u_{Q^{k}}-u_{Q^{j}}|\leq
(k-j)\left(1+2^{n+2}\right)\|u\|_{BMO_{p}}.
\end{equation}
\end{lemma}
\noindent {\bf Proof.} We easily remark that:
$$|Q^{j+1}|=2^{n+2}|Q^{j}|.$$
We compute:
\begin{eqnarray*}
|u_{Q^{j+1}}-u_{Q^{j}}|&=&\frac{1}{|Q^{j}|}\int_{Q^{j}}|u_{Q^{j+1}}-u_{Q^{j}}|\\
&\leq&\frac{1}{|Q^{j}|}\int_{Q^{j}}|u-u_{Q^{j}}|+\frac{1}{|Q^{j}|}\int_{Q^{j}}|u-u_{Q^{j+1}}|\\
&\leq&
\|u\|_{BMO_{p}}+\frac{2^{n+2}}{|Q^{j+1}|}\int_{Q^{j+1}}|u-u_{Q^{j+1}}|\\
&\leq& \|u\|_{BMO_{p}}+2^{n+2}\|u\|_{BMO_{p}}\leq\left(1+2^{n+2}\right)\|u\|_{BMO_{p}},
\end{eqnarray*}
which immediately gives (\ref{meanest_eq1}), and consequently
(\ref{EfS_eq1}). $\hfill{\blacksquare}$\\

\noindent The following two lemmas are of notable importance for the proof of the logarithmic
Sobolev inequality (\ref{cara_eq2}). In the first lemma we bound the
terms  $\phi_{j}*u$ for $j\geq 1$, while, in the second
lemma, we give a bound on $\phi_{0}*u$.
\begin{lemma}{\bf (Estimate of $\|\phi_{j}*u\|_{L^{\infty}(\R^{n+1})}$
    for $j\geq 1$)}\label{key_lemma}\\
Let $u\in BMO_{p}(\R^{n+1})$. Then there exists a constant $C=C(n)>0$
such that:
\begin{equation}\label{key_eq1}
\|u*\phi_{j}\|_{L^{\infty}(\R^{n+1})}\leq C
\|u\|_{BMO_{p}(\R^{n+1})}\quad \mbox{for any}\quad j\geq 1,
\end{equation}
where $(\phi_{j})_{j\geq 1}$ is the sequence of functions given in
(\ref{hedshsu1}).
\end{lemma}
\noindent {\bf Proof.} We will show that
\begin{equation}\label{kersha}
|(\phi_{j}*u)(z)|\leq C \|u\|_{BMO_{p}}\quad \mbox{for}\quad
z=0.
\end{equation}
The general case with $z\in \R^{n+1}$ could be deduced from
(\ref{kersha}) by translation. Throughout the proof, we will sometimes
omit (when there is no confusion) the dependence of the norm on the space $\R^{n+1}$. The
proof is divided into three steps.\\

\noindent {\bf Step 1.} (Decomposition of $(\phi_{j}*u)(0)$ on parabolic
cubes).\\

\noindent Since $\h{\phi}_{j}$ is supported in $\{z\in \R^{n+1};\,
2^{j-1}\leq \|z\|\leq 2^{j+1}\}$ then
$\h{\phi}_{j}(0)=0=\int_{\R^{n+1}}\phi_{j}$. 
Using this equality, we can write:
\begin{equation*}
(\phi_{j}*u)(0)=\intws \phi_{j}(-z)(u(z)-u_{Q^{1-j}}) dz
\end{equation*}
where $Q^{1-j}$ is the parabolic cube defined by (\ref{key_para}) and
centered at $0$. This implies that
\begin{equation}\label{key_eq2}
|(\phi_{j}*u)(0)| \leq  \overbrace{\int_{Q^{1-j}} |\phi_{j}(-z)||u(z)-u_{Q^{1-j}}| dz}^{A_{1}} +
\overbrace{\int_{\R^{n+1}\setminus Q^{1-j}} |\phi_{j}(-z)||u(z)-u_{Q^{1-j}}| dz}^{A_{2}}.
\end{equation}
{\bf Step 1.1.} (Estimate of $A_{1}$).\\

\noindent From (\ref{hedshsu1}), the term $A_{1}$ can be estimated as
follows:
\begin{eqnarray*}
A_{1} &\leq&
2^{(n+2)(j-1)}\|\phi_{1}\|_{L^{\infty}}\int_{Q^{1-j}}|u(z)-u_{Q^{1-j}}|
dz\\
&\leq& 2^{(n+2)(j-1)}|Q^{1-j}|\,\|\phi_{1}\|_{L^{\infty}}\|u\|_{BMO_{p}}\\
&\leq& |Q_{1}|\|\phi_{1}\|_{L^{\infty}}\|u\|_{BMO_{p}},
\end{eqnarray*}
hence
\begin{equation}\label{key_eq-1}
A_{1} \leq C_{0} \|u\|_{BMO_{p}}\quad \mbox{with}\quad
C_{0}=|Q_{1}|\|\phi_{1}\|_{L^{\infty}(\R^{n+1})}.
\end{equation}
{\bf Step 2.} (Estimate of $A_{2}$).\\

\noindent We rewrite $A_{2}$ as the following series:
\begin{equation}\label{key_eq3}
A_{2}=2^{(n+2)(j-1)}\sum_{-\infty<k\leq j} \int_{Q^{2-k}\setminus Q^{1-k}}
\left|\phi_{1}\left(-2^{(j-1)a}z\right)\right| |u(z)-u_{Q^{1-j}}| dz.
\end{equation}
Since $\phi_{1}$ is the inverse Fourier transform of a compactly
supported function then we have:
\begin{equation}\label{key_eq4}
\forall \o{m}\in \N^{*},\;\exists C_{1}>0,\;|\phi_{1}(z)|\leq
\frac{C_{1}}{\|z\|^{\o{m}}} \quad \mbox{for all}\quad \|z\|\geq 1. 
\end{equation}
The asymptotic behavior of $\phi_{1}$ shown by (\ref{key_eq4}) leads to
the following decomposition of the term $A_{2}$:
\begin{eqnarray*}
A_{2}&\leq& \overbrace{C_{1}2^{(n+2)(j-1)} \sum_{-\infty<k\leq
  j}\int_{Q^{2-k}\setminus Q^{1-k}}
\frac{1}{\|2^{(j-1)a}z\|^{\o{m}}}|u(z)-u_{Q^{2-k}}| dz}^{A_{3}}\\
&+& \overbrace{C_{1}2^{(n+2)(j-1)} \sum_{-\infty<k\leq
  j}\int_{Q^{2-k}\setminus
  Q^{1-k}}\frac{1}{\|2^{(j-1)a}z\|^{\o{m}}}|u_{Q^{2-k}}-u_{Q^{1-j}}| dz}^{A_{4}}.
\end{eqnarray*}
\noindent {\bf Step 2.1.} (Estimate of $A_3$).\\

\noindent Since the integral appearing in $A_{3}$ is done over
$Q^{2-k}\setminus Q^{1-k}$, we obtain
$$\|2^{(j-1)a}z\|^{\o{m}}\geq 2^{\o{m}(j-k)}.$$
Using this inequality together with the fact that
$$\int_{Q^{2-k}\setminus Q^{1-k}} |u(z)-u_{Q^{2-k}}| dz\leq
2^{(n+2)(2-k)}|Q_{1}|\|u\|_{BMO_{p}},$$
we can estimate the term $A_{3}$ as follows:
\begin{equation}\label{key_eq8}
A_{3}\leq C_{1}2^{n+2}\left(\sum_{-\infty<k\leq
  j}2^{-(\o{m}-(n+2))(j-k)} \right)|Q_{1}|\|u\|_{BMO_{p}},
\end{equation}
where the above series converges for $\o{m}>n+2$.\\

\noindent {\bf Step 2.2.} (Estimate of $A_4$).\\

\noindent Using Lemma \ref{mean_est}, and the fact that
$\|2^{(j-1)a}z\|^{\o{m}}\geq 2^{\o{m}(j-k)}$ on $Q^{2-k}\setminus
Q^{1-k}$, the term $A_{4}$ can be estimated as follows:
\begin{eqnarray}\label{key_eq9}
A_{4} &\leq& C_{1}2^{(n+2)(j-1)}\left(\sum_{-\infty<k\leq
  j} 2^{-\o{m}(j-k)} (1+j-k) |Q^{2-k}|\right)\|u\|_{BMO_{p}}\nonumber\\
&\leq & C_{1}2^{n+2}\left(\sum_{-\infty<k\leq j}
  2^{-(\o{m}-(n+2))(j-k)}(1+j-k)\right)|Q_{1}|\|u\|_{BMO_{p}},
\end{eqnarray}
where  the above series also converges for $\o{m}>n+2$.\\

\noindent {\bf Step 3.} (Conclusion).\\

\noindent From (\ref{key_eq-1}), (\ref{key_eq8}) and
(\ref{key_eq9}), inequality (\ref{key_eq1}) directly follows with a
constant $C>0$ independent of $j$. $\hfill{\blacksquare}$
\begin{lemma}{\bf (Estimate of
    $\|\phi_{0}*u\|_{L^{\infty}(\R^{n+1})}$)}\label{key_lemma1}\\
Let $u\in W_{2}^{2m,m}(\R^{n+1})$ with $m>\frac{n+2}{4}$. Then there exists a
constant $C=C(m,n)>0$
such that we have:
\begin{equation}\label{lonless}
\|\phi_{0}*u\|_{L^{\infty}}\leq C \left(1+
  \|u\|_{BMO_{p}}\left(1+\log^{+}\|u\|_{W_{2}^{2m,m}}
  \right)\right).
\end{equation}
\end{lemma}
\noindent {\bf Proof.} The constants that will appear may differ from
line to line, but only depend on $n$ and $m$. The
proof of this lemma combines somehow the
proof of Lemmas \ref{BSI_lelemme} and \ref{key_lemma}. We write down
$u_{Q^{1}}$ as a finite sum of a telescopic sequence for $N\geq 1$:
$$ u_{Q^{1}} = (u_{Q^{1}}-u_{Q^{2}})+ \cdots + (u_{Q^{N-1}}-u_{Q^{N}}) +
u_{Q^{N}}.$$
From Lemma \ref{mean_est}, we deduce that:
$$|u_{Q^{1}}|\leq C(N-1)\|u\|_{BMO_{p}}+|u_{Q^{N}}|.$$
Remark that applying Cauchy-Schwarz inequality, we get
$$|u_{Q^{N}}|\leq
\frac{1}{|Q^{N}|}\int_{Q^{N}}|u|\leq
\left(\int_{Q^{N}}u^{2}\right)^{1/2}\left(\int_{Q^{N}}1^{2}\right)^{1/2},$$  
then we obtain
\begin{equation}\label{rf_1}
|u_{Q^{1}}|\leq C\left(N\|u\|_{BMO_{p}} + 2^{-\g N} \|u\|_{W_{2}^{2m,m}}
\right)\quad \mbox{with}\quad \g=\frac{n+2}{2}.
\end{equation}
Following similar arguments as in the proof of Lemma \ref{BSI_lelemme}, we
may optimize (\ref{rf_1}) in $N$, we finally get:
\begin{equation}\label{rf_2}
|u_{Q^{1}}|\leq C \left( 1 + \|u\|_{BMO_{p}}\left( 1 + \log^{+}
    \|u\|_{W_{2}^{2m,m}}\right)\right).
\end{equation}
We now estimate $|(\phi_{0}*u)(z)|$ for $z=0$. Again, the same estimate
could be obtained for any $z\in \R^{n+1}$ by translation. We write
\begin{eqnarray*}
(\phi_{0}*u)(0)&=& \int_{\R^{n+1}} \phi_{0}(-z) u(z)\\
&=& \int_{\R^{n+1}} \phi_{0}(-z) (u(z)-u_{Q^{1}}) + \int_{\R^{n+1}}
\phi_{0}(-z)u_{Q^{1}}\\
&=& \overbrace{\int_{Q_{1}} \phi_{0}(-z) (u(z)-u_{Q^{1}})}^{B_1} +\overbrace{\int_{\R^{n+1}\setminus
  Q^{1}} \phi_{0}(-z) (u(z)-u_{Q^{1}})}^{B_2} + \overbrace{\int_{\R^{n+1}}
\phi_{0}(-z)u_{Q^{1}}}^{B_{3}},
\end{eqnarray*}
where
\begin{equation}\label{wg3a1}
|B_{1}|\leq C \|u\|_{BMO_{p}},
\end{equation}
and, from (\ref{rf_2}),
\begin{equation}\label{wg3a2}
|B_{3}|\leq C  \left( 1 + \|u\|_{BMO_{p}}\left( 1 + \log^{+}
    \|u\|_{W_{2}^{2m,m}}\right)\right).
\end{equation}
In order to estimate $B_{2}$, we argue as in Step 2 of Lemma
\ref{key_lemma}. In fact we have:
\begin{eqnarray}\label{wg3a3}
|B_{2}| &\leq& \sum_{k\geq 1} \int_{Q^{k+1}\setminus Q^{k}}
|\phi_{0}(-z)| |u(z)-u_{Q^{k+1}}| + \sum_{k\geq 1} \int_{Q^{k+1}\setminus Q^{k}}
|\phi_{0}(-z)| |u_{Q^{k+1}}-u_{Q^{1}}|\nonumber\\
&\leq& \left(\sum_{k\geq 1} (\sup_{Q^{k+1}\setminus Q^{k}}|\phi_{0}(-z)|)
|Q^{k+1}| (1+k)\right)\|u\|_{BMO_{p}}\nonumber\\
&\leq& 2^{n+2}\left(\sum_{k\geq 1}
  2^{-(\o{m}-(n+2))}(1+k)\right)|Q_{1}|\|u\|_{BMO_{p}}, 
\end{eqnarray}
where for the last line we have used the fact that $|\phi_{0}(z)|\leq
\frac{C}{\|z\|^{\o{m}}}$ for $\|z\|\geq 1$. Of course the
above series converges if we choose $\o{m}>n+2$. From (\ref{wg3a1}),
(\ref{wg3a2}) and (\ref{wg3a3}), the
result follows. $\hfill{\blacksquare}$

\begin{corollary}{\bf (A control of $\|u\|_{\wt{F}^{0}_{\infty,2}}$)}\label{cara_lem1}\\
Let $u\in BMO_{p}(\R^{n+1})\cap \wt{F}^{0}_{\infty,1}(\R^{n+1})$, then $u\in
\wt{F}^{0}_{\infty,2}(\R^{n+1})$ and we have:
\begin{equation}\label{cara_eq1}
\|u\|_{\wt{F}^{0}_{\infty,2}}\leq C
\|u\|^{1/2}_{BMO_{p}}\|u\|^{1/2}_{\wt{F}^{0}_{\infty,1}},
\end{equation}
where $C=C(n)>0$ is a positive constant.
\end{corollary}
\noindent {\bf Proof.}
Using (\ref{key_eq1}), we compute:
\begin{equation*}
\|u\|_{\wt{F}^{0}_{\infty,2}}=\left\|\left(\sum_{j\geq 1}|\phi_{j}*u|^{2}
  \right)^{1/2} \right\|_{L^{\infty}}\leq \left\|\left(\sup_{j\geq
    1}\|\phi_{j}*u\|_{L^{\infty}}\sum_{j\geq 1}|\phi_{j}*u|
\right)^{1/2}\right\|_{L^{\infty}} \leq C
\|u\|^{1/2}_{BMO_{p}}\|u\|^{1/2}_{\wt{F}^{0}_{\infty,1}},
\end{equation*}
which terminates the proof. $\hfill{\blacksquare}$
\begin{rem}
From \cite{Bownik07}, it seems that $BMO_p$ spaces can be characterized
in terms of parabolic Lizorkin-Triebel spaces. In the case of elliptic
spaces, it is a well-known result (see \cite{Tri92, FraJaw90}) which
allows to simplify the proof of the Kozono-Taniuchi inequality.
\end{rem}

We can now give the proof of our first main result (Theorem \ref{theorem1}).\\

\noindent {\bf Proof of Theorem \ref{theorem1}.} Using (\ref{BSob_ineq})
and (\ref{cara_eq1}), we
obtain:
\begin{equation}\label{cara_eq3}
\|u\|_{\wt{F}^{0}_{\infty,1}}\leq C\left(1+
  \|u\|^{1/2}_{BMO_{p}}\|u\|^{1/2}_{\wt{F}^{0}_{\infty,1}}\left(1+
(\log^{+}\|u\|_{W_{2}^{2m,m}})^{1/2}\right)\right).
\end{equation}
Notice that the constant $C$ can always be chosen such that $C\geq 1$. If
$\|u\|_{\wt{F}^{0}_{\infty,1}}\leq 1$, we evidently have:
\begin{equation}\label{cara_eq4}
\|u\|_{\wt{F}^{0}_{\infty,1}}\leq C\leq C \left(1+ \|u\|_{BMO_{p}} \left(1+
    \log^{+} \|u\|_{W_{2}^{2m,m}}\right) \right).
\end{equation}
If $\|u\|_{\wt{F}^{0}_{\infty,1}}>1$, then, dividing (\ref{cara_eq3}) by
$\|u\|^{1/2}_{\wt{F}^{0}_{\infty,1}}$, we can easily deduce inequality
(\ref{cara_eq4}).  Using the fact that 
$$\|u\|_{L^{\infty}}\leq C\sum_{j\geq 0}\|\phi_{j}*u\|_{L^{\infty}}\leq
C \left( \|\phi_{0}*u\|_{L^{\infty}} + \|u\|_{\wt{F}^{0}_{\infty,1}}\right),$$
and using inequalities (\ref{lonless}) and (\ref{cara_eq4}), we directly
get into the result.$\hfill{\blacksquare}$
\section{A parabolic Kozono-Taniuchi inequality on a bounded domain}\label{sec3}
The goal of this section is to present, on the one hand, the proof of Theorem
\ref{theorem2}. On the other hand, at the end of this section, we give an application where
we show how to use inequality (\ref{theo2_eq1}) in order to maintain the
long-time existence of solutions to some parabolic equations. Let us indicate that
throughout this section, the positive constant $C=C(T)>0$ may vary from
line to line.
\subsection{Proof of Theorem \ref{theorem2}}
In order to simplify the arguments of the proof, we first show Theorem
\ref{theorem2} in the special case when $n=m=1$. Then
we give the principal ideas how to prove the result in the general
case. Call
$$I=(0,1)\quad \mbox{and}\quad \Omega_{T}=I\times (0,T),$$
we first show the following proposition:
\begin{proposition}{\bf (Theorem \ref{theorem2}, case: $n=m=1$)}\label{prop1}\\
Let $u\in W_{2}^{2,1}(\Omega_{T})$. Then there exists
a constant $C=C(T)>0$ such that:
\begin{equation}\label{prop1_eq}
\|u\|_{L^{\infty}(\Omega_{T})}\leq C \left(1+
  \|u\|_{\overline{BMO}_{p}(\Omega_{T})} \left(1+
    \log^{+} \|u\|_{W_{2}^{2,1}(\Omega_{T})}\right) \right).
\end{equation}
\end{proposition}
As a similar inequality of (\ref{prop1_eq}) is already shown on
$\R^{2}$ (see inequality (\ref{cara_eq2})), the idea of
the proof of (\ref{prop1_eq}) lies in using (\ref{cara_eq2}) for a
special extension of the function $u\in W_{2}^{2,1}(\Omega_{T})$ to the entire
space $\R^{2}$. For this reason, we demand that the extended function stays in
$W_{2}^{2,1}(\R^{2})$ which is done via the following arguments. Remark
first that the function $u$ can be extended by continuity to the
boundary $\partial \Omega_{T}$ of $\Omega_{T}$. Take $\tilde{u}$ as the function defined over
$$\widetilde{\Omega}_{T}=(-1,2)\times (-T,2T)$$
as follows:
\begin{equation}\label{first_ext}
\tilde{u}(x,t)=
\left\{
\begin{aligned}
& -3u(-x,t)+4u\left(-\frac{x}{2},t\right)\quad&\mbox{for}&\quad &-1<x<0,\, 0\leq
t\leq T,&\\
& -3 u(2-x,t)+4 u\left(\frac{3-x}{2},t\right)\quad&\mbox{for}&\quad &1<x<2,\, 0\leq
t\leq T,&
\end{aligned}
\right.
\end{equation}
and
\begin{equation}\label{second_ext}
\tilde{u}(x,t)=
\left\{
\begin{aligned}
& u(x,-t)\quad &\mbox{for}&\quad &-T<t\leq 0&\\
& u(x,2T-t)\quad &\mbox{for}&\quad &T\leq t< 2T&.
\end{aligned}
\right.
\end{equation}
A direct consequence of this extension is the following
lemma.
\begin{lemma}{\bf ($L^{1}$ estimate of $\tilde{u}$)}\label{Hgibi-1}\\
Let $\tilde{u}$ be the function defined by (\ref{first_ext}) and
(\ref{second_ext}). Then there exists a constant $C=C(T)>0$ such
that:
\begin{equation}\label{sa5sD6}
\|\tilde{u}\|_{L^{1}(\widetilde{\Omega}_{T})}\leq C
\|u\|_{L^{1}(\Omega_{T})}.
\end{equation}
\end{lemma}
\noindent {\bf Proof.} The proof of this lemma is direct by the
extension. $\hfill{\blacksquare}$\\

\noindent Another important consequence of the extension (\ref{first_ext}) and
(\ref{second_ext}) is the fact that $\tilde{u}\in
W_{2}^{2,1}(\widetilde{\Omega}_{T})$, and that we have (see for instance
\cite{EVANS})
\begin{equation}\label{5s5}
\|\tilde{u}\|_{W_{2}^{2,1}(\widetilde{\Omega}_{T})}\leq C \|u\|_{W_{2}^{2,1}(\Omega_{T})},\quad C=C(T)>0.
\end{equation}
Let $\mathcal{Z}_{1}\subset\mathcal{Z}_{2}$ be the two subsets of
$\widetilde{\Omega}_{T}$ defined by:
$$\mathcal{Z}_{1}=\{(x,t);\, -1/4<x<5/4\quad\mbox{and}\quad
-T/4<t<5T/4\},$$
and
$$\mathcal{Z}_{2}=\{(x,t);\, -3/4<x<7/4\quad\mbox{and}\quad
-3T/4<t<7T/4\}.$$
Taking the cut-off function $\Psi\in C^{\infty}_{0}(\R^{2})$, $0\leq
\Psi\leq 1$ satisfying:
\begin{equation}\label{cut_off}
\Psi(x,t)=
\left\{
\begin{aligned}
& 1\quad &\mbox{for}&\quad (x,t)\in \mathcal{Z}_{1}\\
& 0\quad &\mbox{for}& \quad (x,t)\in \R^{2}\setminus \mathcal{Z}_{2},
\end{aligned}
\right.
\end{equation}
we can easily deduce from (\ref{5s5}) that $\Psi \tilde{u}\in
W_{2}^{2,1}(\R^{2})$, and
\begin{equation}\label{5s5_1}
\|\Psi\tilde{u}\|_{W_{2}^{2,1}(\R^{2})}\leq C \|u\|_{W_{2}^{2,1}(\Omega_{T})}.
\end{equation}
Since $\Psi \tilde{u}\in W_{2}^{2,1}(\R^{2})$, we can apply inequality
(\ref{cara_eq2}) to the function $\Psi \tilde{u}$, and, having
(\ref{5s5_1}) in hands, the proof of Proposition \ref{prop1} directly
follows if we can show that
\begin{equation}\label{5s5_2}
\|\Psi \tilde{u}\|_{BMO_{p}(\R^{2})}\leq C \|u\|_{\overline{BMO}_{p}(\Omega_{T})},
\end{equation}
and this will be done in the forthcoming arguments.
\subsubsection{Proof of Proposition \ref{prop1}}
In all what follows, it will be useful to deal with an equivalent norm
of the $BMO_p$ space. This norm is given by the following lemma.
\begin{lemma}{\bf (Equivalent $BMO_{p}$ norms)}\label{equiv_bmo}\\
Let $u\in BMO_{p}(\mathcal{O})$, $\mathcal{O}\subseteq \R^{n+1}$ is an
open set. The parabolic $BMO_p$ norm of $u$ given by (\ref{bmo_norm}) is equivalent
to the following norm for which we give the same notation:
\begin{equation}\label{eq_bmo_nor}
\|u\|_{BMO_{p}(\mathcal{O})}=\sup_{Q\subset \mathcal{O}} \left(\inf_{c\in \R}
\frac{1}{|Q|} \int_{Q}|u-c|\right),\quad Q\mbox{ given by (\ref{haAnas})}.
\end{equation}
\end{lemma}
{\bf Proof.} The proof of this lemma is direct. It suffices to see that
for any $c\in \R$, we have:
$$|u-u_{Q}|\leq |u-c|+|c-u_{Q}|\leq |u-c|+ \frac{1}{|Q|}
\int_{Q}|u-c|,$$
which immediately gives:
$$\int_{Q}|u-u_{Q}|\leq 2 \int_{Q}|u-c|,$$
hence
\begin{equation}\label{Y5YL}
\frac{1}{2|Q|}\int_{Q}|u-u_{Q}|\leq\inf_{c\in
  \R}\frac{1}{|Q|}\int_{Q}|u-c|\leq \frac{1}{|Q|}\int_{Q}|u-u_{Q}|,
\end{equation}
and the equivalence of the two norms follows. $\hfill{\blacksquare}$\\

\noindent From now on, and for the sake of simplicity, we will denote:
$$\inm_{Q}u=\frac{1}{|Q|}\int_{Q}u .$$
The following lemma gives an estimate of $\displaystyle \inf_{c\in
  \R}\inm_{Q}|u-c|$ on small parabolic cubes.
\begin{lemma}\label{Hgibi1}
Let $f\in L^{1}_{loc}(\R^{2})$. Take $Q_{r}\subseteq Q_{2r}$ two
parabolic cubes of $\R^{2}$. We do
not require that the cubes have the same center. Then we have:
\begin{equation}\label{moa1}
\inf_{c\in \R}\inm_{Q_{r}}|f-c|\leq 8\inf_{c\in\R}\inm_{Q_{2r}}|f-c|.
\end{equation}
\end{lemma}
{\bf Proof.} For $c\in \R$, we compute:
$$\inm_{Q_{r}}|f-c|\leq \frac{|Q_{2r}|}{|Q_{r}|} \inm_{Q_{2r}}|f-c|\leq 8
\inm_{Q_{2r}}|f-c|.$$
Taking the infimum of both sides we arrive to the
result. $\hfill{\blacksquare}$\\

\noindent The next lemma gives an estimate of $\displaystyle
\inf_{c\in\R}\inm_{Q_{r}}|\tilde{u}-c|$ on small parabolic cubes in
$$\widehat{\Omega}_{T}=(-1,2)\times (0,T).$$
Define the term $r_{0}>0$ as the greatest positive real number such that
there exists $Q_{r_{0}}\subseteq \Omega_{T}$, i.e.,
\begin{equation}\label{cond_on_r}
r_{0}=\sup\{r>0;\, r\leq 1/2 \mbox{ and } r^{2}\leq T/2\}.
\end{equation}
We show the following:
\begin{lemma}{\bf (Estimates on small parabolic cubes in
    $\widehat{\Omega}_{T}$)}\label{Hgibi2}\\
Let $\tilde{u}$ be the function defined by (\ref{first_ext}) and
(\ref{second_ext}). Take any parabolic cube $Q_{r}$ satisfying:
\begin{equation}\label{moa2}
Q_{r}\subseteq \widehat{\Omega}_{T},\quad \mbox{with}
\quad r\leq r_{1}\mbox{ and }2r_{1}=r_{0},
\end{equation}
where $r_{0}$ is given by (\ref{cond_on_r}). Then there exists a
universal constant $C>0$ such that:
\begin{equation}\label{moa3}
\inf_{c\in\R}\inm_{Q_{r}}|\tilde{u}-c|\leq C\|u\|_{BMO_{p}(\Omega_{T})}.
\end{equation}
\end{lemma}
{\bf Proof.} Call $\Omega^{d}_{T}$ and $\Omega^{g}_{T}$ the right and
the left neighbor sets of $\Omega_T$ defined respectively by:
$$\Omega^{d}_{T}=(-1,0)\times (0,T)\quad \mbox{and} \quad
\Omega^{g}_{T}=(1,2)\times (0,T).$$
First let us mention that if the cube $Q_{r}$ lies in
$\Omega_{T}$ then inequality (\ref{moa3}) is evident (see the equivalent
definition (\ref{eq_bmo_nor}) of the parabolic $BMO_{p}$ norm). Two
remaining cases are to be considered: either $Q_{r}$ intersects the set
$\{x=0\}\cup\{x=1\}$, or $Q_{r}$ lies in $\Omega^{d}_{T}\cup
\Omega^{g}_{T}$. Our assumption (\ref{moa2}) on the radius of the
parabolic cube makes it impossible that the cube $Q_{r}$ meets
$\Omega^{d}_{T}$ and $\Omega^{g}_{T}$ at the same time. Therefore, and
in order to make the proof simpler, we only consider the following
cases: either $Q_{r}$ intersects the set $\{x=0\}$, or $Q_{r}$ lies in
$\Omega^{g}_{T}$. The proof is then divided
into three main steps:\\

\noindent {\bf Step 1.} ($Q_{r}$ intersects the line $\{x=0\}$).\\

\noindent {\bf Step 1.1.} (First estimate).\\

\begin{figure}[!h]
\psfrag{-1}{\footnotesize{$-1$}}
\psfrag{0}{\footnotesize{$0$}}
\psfrag{OT}{\footnotesize{$\Omega_T$}}
\psfrag{OTr}{\footnotesize{$\Omega^{d}_{T}$}}
\psfrag{OTl}{\footnotesize{$\Omega^{g}_{T}$}}
\psfrag{1}{\footnotesize{$1$}}
\psfrag{T}{\footnotesize{$T$}}
\psfrag{2}{\footnotesize{$2$}}
\psfrag{x}{\footnotesize{$x$}}
\psfrag{t}{\footnotesize{$t$}}
\psfrag{Qr}{\footnotesize{$Q_{r}$}}
\psfrag{Q2r}{\footnotesize{$\!\!Q_{2r}$}}
\centering\epsfig{file=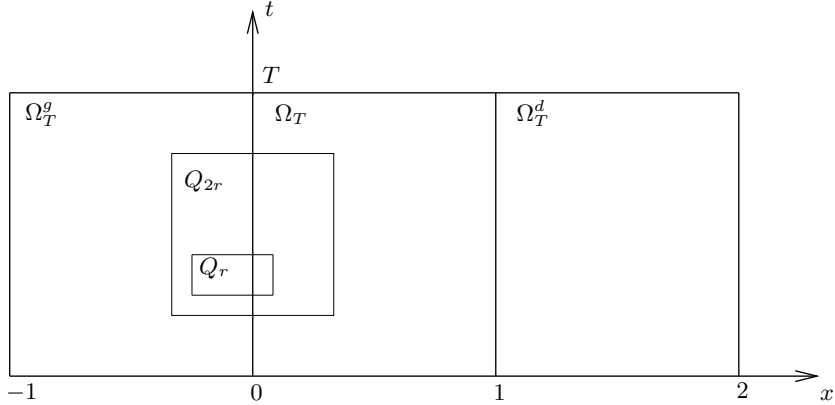, width=110mm}
   \caption{Analysis on cubes intersecting $\{x=0\}$}\label{cubes}
\end{figure}

\noindent Again the assumption (\ref{moa2}) imposed on the radius $r$
makes it possible to embed $Q_{r}$ in a larger parabolic cube
$Q_{2r}\subseteq \widehat{\Omega}_{T}$ of
radius $2r$, which is symmetric with respect to the line $\{x=0\}$ (see
Figure \ref{cubes}).
Then the center of the cube $Q_{2r}$ should be also on the same
line, but we do not require that the two cubes
$Q_{r}$ and $Q_{2r}$ have centers with the same ordinate $t$. Now, using
Lemma
\ref{Hgibi1}, we deduce that:
\begin{equation}\label{cri1}
\inf_{c\in\R}\inm_{Q_{r}}|\tilde{u}-c|\leq 8
\inf_{c\in\R}\inm_{Q_{2r}}|\tilde{u}-c|,
\end{equation}
and hence in order to conclude, we need to estimate the right-hand side
of the above inequality with respect to
$\|u\|_{BMO_{p}(\Omega_{T})}$. Call $Q^{d}_{2r}$ and $Q^{g}_{2r}$ the right and
the left sides of $Q_{2r}$ defined respectively by:
$$Q^{d}_{2r}=Q_{2r}\cap \Omega_{T}\quad \mbox{and} \quad
Q^{g}_{2r}=Q_{2r}\cap \Omega^{g}_{T}.$$
Also call $Q^{{trans}}_{2r}\subseteq \Omega_{T}$, the translation of the cube $Q_{2r}$
by the vector $(2r,0)$, i.e.
$$Q^{{trans}}_{2r}=(2r,0)+ Q_{2r}.$$
For $c\in\R$, we compute:
\begin{eqnarray}\label{moa4}
\int_{Q_{2r}}|\tilde{u}-c|&=&\int_{Q^{g}_{2r}}|\tilde{u}-c|
+\int_{Q^{d}_{2r}}|u-c|,\nonumber\\
&\leq&\int_{Q^{g}_{2r}}|\tilde{u}-c|
+\int_{Q^{trans}_{2r}}|u-c|,
\end{eqnarray}
where we have used the fact that $\tilde{u}=u$ on $\Omega_{T}$, and that
$Q^{d}_{2r}\subseteq Q^{trans}_{2r}$.\\

\noindent {\bf Step 1.2.} (Estimate of
$\int_{Q^{g}_{2r}}|\tilde{u}-c|$).\\

\noindent We compute (using the definition (\ref{first_ext}) of the
function $\tilde{u}$ on $\Omega^{g}_{T}$):
\begin{eqnarray}\label{moa5}
\int_{Q^{g}_{2r}}|\tilde{u}(x,t)-c|dxdt&=&\int_{Q^{g}_{2r}}|-3u(-x,t)+4u(-x/2,t)-c|dxdt\nonumber\\
&\leq& 3\int_{Q^{g}_{2r}}|u(-x,t)-c|dxdt + 4
\int_{Q^{g}_{2r}}|u(-x/2,t)-c|dxdt\nonumber\\
&\leq& 3\int_{Q^{d}_{2r}}|u(x,t)-c|dxdt + 8\int_{Q^{\bar{d}}_{2r}}|u(x,t)-c|dxdt,
\end{eqnarray}
where
$$Q^{\bar{d}}_{2r}=\{(x/2,t);\, (x,t)\in Q^{d}_{2r}\}\subseteq
Q^{d}_{2r}\subseteq Q^{trans}_{2r}.$$
From (\ref{moa5}) we easily deduce that:
$$
\int_{Q^{g}_{2r}}|\tilde{u}-c|\leq 11 \int_{Q^{trans}_{2r}}|u-c|,
$$
and hence (using (\ref{moa4})), we finally get:
\begin{equation}\label{moa6}
\int_{Q_{2r}}|\tilde{u}-c|\leq 12 \int_{Q^{trans}_{2r}}|u-c|.
\end{equation}
Since $|Q_{2r}|=|Q^{trans}_{2r}|$, inequality (\ref{moa6}) gives
$$\inm_{Q_{2r}}|\tilde{u}-c|\leq 12 \inm_{Q^{trans}_{2r}}|u-c|.$$
Since $Q^{trans}_{2r}$ is a parabolic cube in $\Omega_{T}$, taking the
infimum over $c\in\R$ of the above inequality, we obtain:
\begin{equation}\label{moa7}
\inf_{c\in \R}\inm_{Q_{2r}}|\tilde{u}-c|\leq 12
\|u\|_{BMO_{p}(\Omega_{T})}.
\end{equation}
From (\ref{cri1}) and (\ref{moa7}), we deduce (\ref{moa3}).\\

\noindent {\bf Step 2.} ($Q_{r}\subseteq \Omega^{g}_{T}$).\\

\begin{figure}[!h]
\psfrag{-1}{\footnotesize{$\!\!\!\!-1$}}
\psfrag{0}{\footnotesize{$0$}}
\psfrag{OT}{\footnotesize{$\Omega_T$}}
\psfrag{OTr}{\footnotesize{$\Omega^{d}_{T}$}}
\psfrag{OTl}{\footnotesize{$\Omega^{g}_{T}$}}
\psfrag{1}{\footnotesize{$1$}}
\psfrag{a}{\footnotesize{$a_0$}}
\psfrag{b}{\footnotesize{$b_0$}}
\psfrag{c}{\footnotesize{$\!\!\!\!a_1$}}
\psfrag{d}{\footnotesize{$\!\!\!b_1$}}
\psfrag{s}{\footnotesize{$\!\!\!\frac{a_0}{2}$}}
\psfrag{r}{\footnotesize{$\!\!\!\frac{b_0}{2}$}}
\psfrag{-a}{\footnotesize{$\!\!\!\!-a_0$}}
\psfrag{-b}{\footnotesize{$\!\!\!-b_0$}}
\psfrag{T}{\footnotesize{$T$}}
\psfrag{2}{\footnotesize{$2$}}
\psfrag{x}{\footnotesize{$x$}}
\psfrag{t}{\footnotesize{$t$}}
\psfrag{Qr}{\footnotesize{$\!\!Q_{r}$}}
\psfrag{u}{\footnotesize{$\!\!\!Q^{\bar{s}}_{r}$}}
\psfrag{Qrs}{\footnotesize{$\!\!Q^{s}_{r}$}}
\centering\epsfig{file=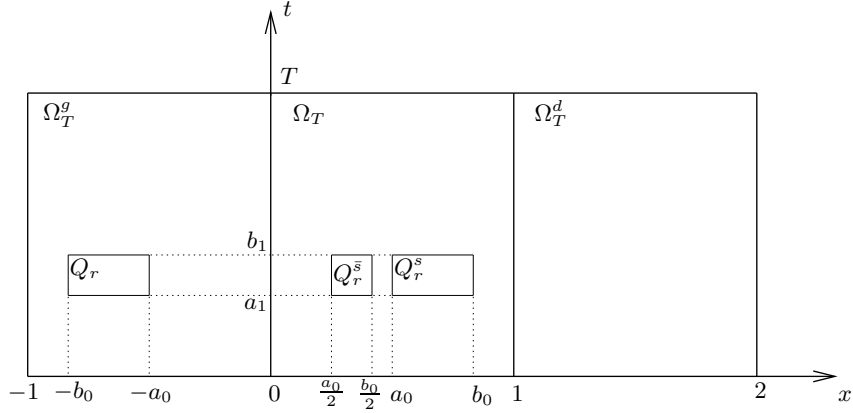, width=110mm}
   \caption{Analysis on cubes $Q_{r}\subseteq \Omega_{T}^{g}$}\label{cubes1}
\end{figure}

\noindent Let $0<a_{0}<b_{0}<1$ and $0<a_{1}<b_{1}<T$ be such that
$$Q_{r}=(-b_{0},-a_{0})\times (a_{1},b_{1}).$$
For any $c\in \R$, we compute:
\begin{eqnarray}\label{mar_g1}
\int_{Q_{r}}|\tilde{u}(x,t)-c|dxdt&=&\int_{Q_{r}}|-3u(-x,t)+4u(-x/2,t)-c|dxdt\nonumber\\
&\leq& 3\int_{Q^{s}_{r}}|u(x,t)-c|dxdt + 8\int_{Q^{\bar{s}}_{r}}|u(x,t)-c|dxdt
\end{eqnarray}
with (see Figure \ref{cubes1}),
$$Q^{s}_{r}=(a_{0},b_{0})\times (a_{1},b_{1})\subseteq \Omega_{T}\quad
\mbox{and}\quad Q^{\bar{s}}_{r}=\left(\frac{a_{0}}{2},\frac{b_{0}}{2}\right)\times
(a_{1},b_{1})\subseteq \Omega_{T}.$$
We remark that $Q^{s}_{r}$ is a parabolic cube in $\Omega_T$, while
$Q^{\bar{s}}_{r}$ is not (its aspect ratio is different). In fact
$Q^{\bar{s}}_{r}$ could be embedded in
a parabolic cube $Q^{\bar{s}}_{r}\subseteq Q^{\bar{\bar{s}}}_{r}\subseteq
\Omega_{T}$, where $Q^{\bar{\bar{s}}}_{r}$ is simply a space translation
of $Q^{s}_{r}$. In particular we have:
\begin{equation}\label{janan7}
|Q_{r}|=|Q^{s}_{r}|=|Q^{\bar{\bar{s}}}_{r}|.
\end{equation}
The above arguments, together with (\ref{mar_g1}) give:
\begin{equation}\label{nhsd}
\int_{Q_{r}}|\tilde{u}-c|\leq 3\int_{Q^{s}_{r}}|u-c| +
8\int_{Q^{\bar{\bar{s}}}_{r}}|u-c|.
\end{equation}
Taking the infimum in $c\in\R$ for both sides of inequality
(\ref{nhsd}), leads to
\begin{equation}\label{janan7bi}
\inf_{c\in\R} \inm_{Q_{r}}|\tilde{u}-c|\leq 11 \|u\|_{BMO_{p}(\Omega_{T})},
\end{equation}
which implies (\ref{moa3}).\\

\noindent {\bf Step 3.} (Conclusion).\\

\noindent As it was already mentioned at the beginning of the proof, the
case where the parabolic cube $Q_{r}$ meets the line $\{x=1\}$ or lies
completely in $\Omega^{d}_{T}$, could be treated using identical
arguments. Therefore, for all small parabolic cubes $Q_{r}$ satisfying
(\ref{moa2}), inequality (\ref{moa3}) is always valid, and this
terminates the proof of Lemma \ref{Hgibi2}. $\hfill{\blacksquare}$\\

\begin{figure}[!h]
     \begin{center}
     \psfrag{-1}{\footnotesize{$\!\!\!\!-1$}}
     \psfrag{0}{\footnotesize{$\!0$}}
     \psfrag{T}{\footnotesize{$\!\!T$}}
     \psfrag{2t}{\footnotesize{$2T$}}
     \psfrag{t}{\footnotesize{$t$}}
     \psfrag{Qr}{\footnotesize{$\!\!Q_{r}$}}
     \psfrag{Qrsy}{\footnotesize{$\!\!Q^{sym}_{r}$}}
     \psfrag{Qrr2}{\footnotesize{$\!\!Q_{r\sqrt{2}}$}}
     \psfrag{2}{\footnotesize{$\!\!2$}}
     \begin{minipage}[b]{.35\linewidth}
      \centering\epsfig{figure=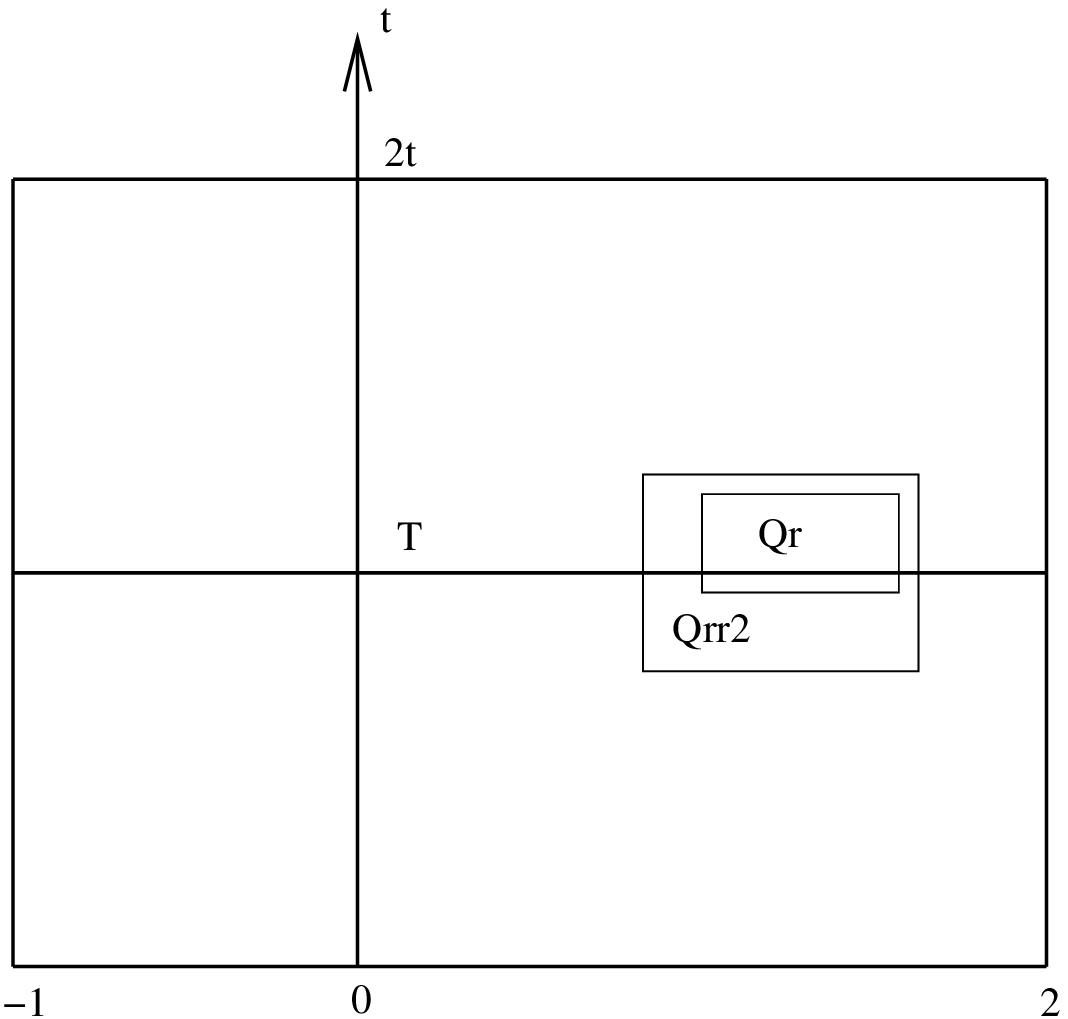,width=\linewidth}
      \caption{$Q_{r}\cap \{t=T\}\neq\emptyset$}\label{cubes_t1}
     \end{minipage}\hspace{1.3cm}
     \begin{minipage}[b]{.35\linewidth}
      \centering\epsfig{figure=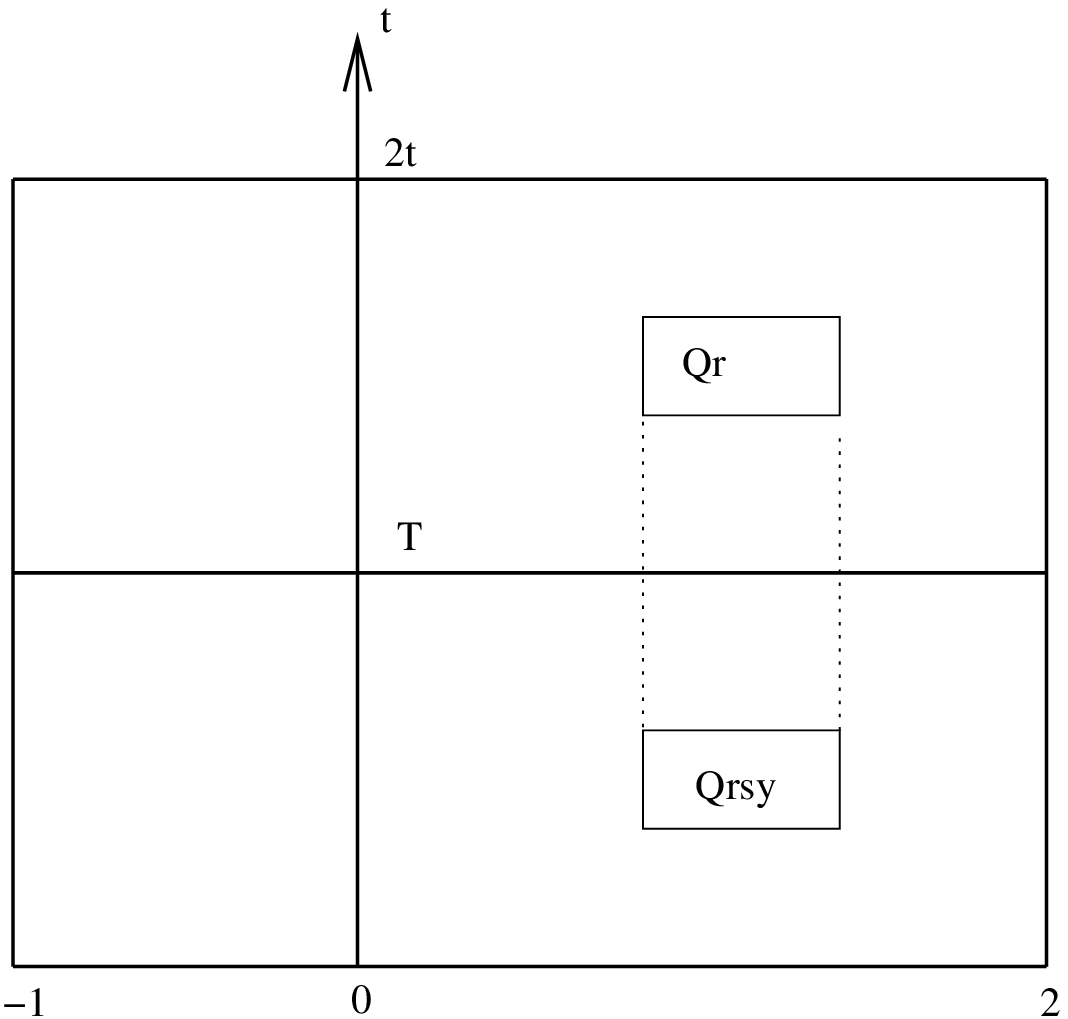 ,width=\linewidth}
      \caption{$Q_{r}\cap \{t=T\}=\emptyset$}\label{cubes_t}
     \end{minipage}
  \end{center}
    \end{figure}

\noindent A generalization of Lemma \ref{Hgibi2} is now given.
\begin{lemma}\label{Hgibi3}{\bf (Estimates on small parabolic cubes in
    $\widetilde{\Omega}_{T}$)}\\
Let $\tilde{u}$ be the function defined by (\ref{first_ext}) and
(\ref{second_ext}). Take any parabolic cube $Q_{r}\subseteq
\widetilde{\Omega}_{T}$ satisfying:
\begin{equation}\label{mabaa0}
r\leq r_{2}\quad \mbox{with} \quad r_{2}\sqrt{2}=r_{1},
\end{equation}
where $r_{1}$ is given by (\ref{moa2}). Then there exists a universal
constant $C>0$ such that:
\begin{equation}\label{mabaa1}
\inf_{c\in\R}\inm_{Q_{r}}|\tilde{u}-c|\leq C\|u\|_{BMO_{p}(\Omega_{T})}.
\end{equation}
\end{lemma}
{\bf Sketch of the proof.} The arguments leading to the proof of this lemma are
already contained in the proof of Lemma \ref{Hgibi2}. First notice that
if $Q_{r}\subseteq \widehat{\Omega}_{T}$, we enter directly (since
$r\leq r\sqrt{2}\leq r_{1}$) to the framework of Lemma \ref{Hgibi2}, and hence
(\ref{mabaa1}) is direct. Because $r\leq r_{1}$, remark that there
exists a cube $Q^{'}_{r}$ obtained by a time translation of $Q_{r}$ such
that $Q^{'}_{r}\subseteq \widehat{\Omega}_{T}$. Therefore it is impossible
that $Q_{r}$ meets at the same time $(-1,2)\times (T,2T)$ and
$(-1,2)\times (-T,0)$. For this reason, we either consider parabolic
cubes intersecting $\{t=T\}$ (see Figure \ref{cubes_t1}),
or parabolic cubes in $(-1,2)\times (T,2T)$ (see Figure \ref{cubes_t}).\\

\noindent \textit{Case $Q_{r}\cap \{t=T\}\neq\emptyset$}. In this case, we
first embed $Q_{r}$ in a larger parabolic cube $Q_{r\sqrt{2}}$ which
is symmetric with respect to the line $\{t=T\}$, so the center of this
cube lies in $\{t=T\}$. We now repeat the same arguments as in Step
1 of Lemma \ref{Hgibi2}, using in particular the symmetry
(\ref{second_ext}) of the function $\tilde{u}$ with respect to
$\{t=T\}$, and the fact that we can consider the cube
$$Q^{trans{'}}_{r\sqrt{2}}=(0,-2r^{2})+Q_{r\sqrt{2}}$$
such that
$$Q^{trans{'}}_{r\sqrt{2}}\subseteq Q_{r_{1}} \subseteq
\widehat{\Omega}_{T},$$
for some cube $Q_{r_{1}}$. Indeed, estimates on all such cubes
$Q^{trans{'}}_{r\sqrt{2}}$ are already controlled by (\ref{moa3}).\\

\noindent \textit{Case $Q_{r}\cap \{t=T\}=\emptyset$}. In this case we
repeat the same arguments as in Step 2 of Lemma \ref{Hgibi2}. Indeed, in
the present case, it is even simpler since the function $\tilde{u}$ is
symmetric with respect to $\{t=T\}$. $\hfill{\blacksquare}$\\

\noindent We now show how to prove estimate (\ref{5s5_2}).\\

\noindent {\bf Proof of estimate (\ref{5s5_2}).} The parabolic $BMO_p$
norm (\ref{eq_bmo_nor}) of $\Psi\tilde{u}$ could be estimated taking the
supremum of
$\inm_{Q_{r}}|\Psi\tilde{u}-(\Psi\ti{u})_{Q_{r}}|$,
$Q_{r}\subseteq \R^{2}$, over
small parabolic cubes ($Q_{r}$ with $r\leq r_{2}/2$), and
big parabolic cubes ($Q_{r}$ with $r>r_{2}/2$). The proof is then divided
into two steps.\\

\noindent {\bf Step 1.} (Analysis on big parabolic cubes $Q_{r}$, $r>r_{2}/2$).\\

\noindent We compute, using the fact that  $\Psi=0$ on
$\R^{2}\setminus \mathcal{Z}_{2}$, and $\Psi\leq
1$ on $\R^{2}$ (see (\ref{cut_off})):
\begin{eqnarray}\label{misang}
\inm_{Q_{r}}|\Psi\tilde{u}-(\Psi\ti{u})_{Q_{r}}|&\leq&
2\inm_{Q_{r}}|\Psi\ti{u}|\leq
\frac{2}{|Q_{r}|}\int_{Q_{r}\cap \mathcal{Z}_{2}}|\ti{u}|\nonumber\\
&\leq& \frac{2^2}{r^{3}_{2}}\int_{Q_{r}\cap \mathcal{Z}_{2}}|\ti{u}|\leq
\frac{2^2}{r^{3}_{2}}\int_{\widetilde{\Omega}_{T}}|\ti{u}|\leq
C\|u\|_{L^{1}(\Omega_{T})}.
\end{eqnarray}
\noindent {\bf Step 2.} (Analysis on small parabolic cubes $Q_{r}$,
$r\leq r_{2}/2$).\\

\noindent From the definition (\ref{mabaa0}) of $r_{2}$, and the
construction (\ref{cut_off}) of the function $\Psi$, we deduce that if
$Q_{r}$ intersects $\mathcal{Z}_{2}$ then forcedly $Q_{r}\subseteq
\widetilde{\Omega}_{T}$. If not, i.e. $Q_{r}\cap \mathcal{Z}_{2}=\emptyset$
then $\Psi=0$ on $Q_{r}$, and therefore:
\begin{equation}\label{result2}
\inm_{Q_{r}}|\Psi\tilde{u}-(\Psi\ti{u})_{Q_{r}}|=0.
\end{equation}
Then we have only to  consider $Q_{r}\subseteq
\widetilde{\Omega}_{T}$.\\ 

\noindent {\bf Step 2.1.} (First estimate).\\

\noindent Using (\ref{Y5YL}), we get
\begin{equation}\label{Reg_eq1}
\inm_{Q_{r}}|\Psi\ti{u} - (\Psi\ti{u})_{Q_{r}}|\leq 2\inf_{c\in\R}\inm_{Q_{r}}|\Psi\ti{u}-c|\leq
2\inm_{Q_{r}}|\Psi\ti{u}-c_{0}\Psi_{Q_{r}}|, 
\end{equation}
for any fixed constant $c_{0}\in \R$. Remark that we can write:
\begin{equation}\label{Reg_eq2}
\Psi\ti{u}-c_{0}\Psi_{Q_{r}}=(\Psi-\Psi_{Q_{r}})\ti{u} +
 (\ti{u}-c_{0})\Psi_{Q_{r}}.
\end{equation}
Hence, we deduce that
\begin{eqnarray}\label{Reg_eq3}
\inm_{Q_{r}}|\Psi\ti{u} - (\Psi\ti{u})_{Q_{r}}|&\leq& Cr\inm_{Q_{r}}|\ti{u}|+
2\inf_{c_{0}\in\R}\inm_{Q_{r}}|\ti{u}-c_{0}|\nonumber\\ 
&\leq& Cr\inm_{Q_{r}}|\ti{u}|+ 2C \|u\|_{BMO_{p}(\Omega_{T})},
\end{eqnarray}
where for the first line we have used that fact that $\Psi \leq 1$ and
that $\Psi$ is Lipschitz, and for the second line we have used
(\ref{mabaa1}).\\

\noindent {\bf Step 2.2.} (Estimate of $\inm_{Q_{r}}|\ti{u}|$).\\ 

\noindent We have 
\begin{eqnarray}\label{Azmeteq}
\inm_{Q_{r}}|\ti{u}|&\leq& |\ti{u}_{Q_{r}}| + \inm_{Q_{r}} |\ti{u} -
\ti{u}_{Q_{r}}|\nonumber\\
&\leq& |\ti{u}_{Q_{r}}| + 2 \inf_{c\in\R}\inm_{Q_{r}}|\ti{u}-c|\nonumber\\
&\leq& |\ti{u}_{Q_{r}}| + 2C\|u\|_{BMO_{p}(\Omega_{T})},
\end{eqnarray}
where for the second line, we have used (\ref{Y5YL}), while for the
third line, we have used (\ref{mabaa1}). Remark that from the proof of
Lemma \ref{mean_est} with $n=1$, we have for $Q_{2^{j}r}\subseteq
Q_{2^{j+1}r}\subseteq \widetilde{\Omega}_{T}$:
\begin{eqnarray*}
|\ti{u}_{Q_{2^{j}r}}-\ti{u}_{Q_{2^{j+1}r}}|&\leq&
\inm_{Q_{2^{j}r}}|\ti{u} - \ti{u}_{Q_{2^{j}r}}| + 2^{3}
\inm_{Q_{2^{j+1}r}}|\ti{u}-\ti{u}_{Q_{2^{j+1}r}}|\\ 
&\leq& 2(1+2^{3})\sup_{Q_{\rho}\subseteq
  \widetilde{\Omega}_{T},\,\rho\leq 2^{j+1}r}
\left(\inf_{c\in\R}\inm_{Q_{\rho}}|\ti{u}-c|\right)\\
&\leq& 2C(1+2^{3}) \|u\|_{BMO_{p}(\Omega_{T})},
\end{eqnarray*}
where we have used (\ref{Y5YL}) for the second line, and, for the third
line, we have used (\ref{mabaa1}) assuming $2^{j+1}r\leq
r_{2}$. Defining
$$j_{0}=\min\{j\in \N;\, r_{2}/2\leq 2^{j}r < r_{2}\},$$ 
and using a telescopic sequence, we can deduce that 
\begin{eqnarray}\label{Reg_eq4}
|\ti{u}_{Q_{r}} - \ti{u}_{Q_{2^{j_{0}}r}}| &\leq& j_{0} (2C(1+2^{3}))
\|u\|_{BMO_{p}(\Omega_{T})}\nonumber\\
&\leq& C(1+|\log r|)\|u\|_{BMO_{p}(\Omega_{T})}. 
\end{eqnarray}
Moreover, we have
\begin{equation}\label{JM1}
|\ti{u}_{Q_{2^{j_{0}}r}}|\leq
\frac{1}{|Q_{r_{2}/2}|}\int_{\widetilde{\Omega}_{T}} |\ti{u}|\leq
C\|u\|_{L^{1}(\Omega_{T})}, 
\end{equation} 
where we have used (\ref{sa5sD6}) for the second inequality. From
(\ref{Reg_eq3}), (\ref{Reg_eq4}) and (\ref{JM1}), we get:
\begin{equation}\label{JM2}
\inm_{Q_{r}}|\ti{u}|\leq C \left(\|u\|_{L^{1}(\Omega_{T})} + (1+|\log
  r|) \|u\|_{BMO_{p}(\Omega_{T})} \right)
\end{equation}
for some constant $C>0$.\\

\noindent {\bf Step 2.3.} (Conclusion for $r\leq r_{2}/2$).\\

\noindent Finally, putting together (\ref{Reg_eq3}) and (\ref{JM2}), we
deduce that 
\begin{eqnarray}\label{mwbasyd}
\inm_{Q_{r}}|\Psi\ti{u} - (\Psi\ti{u})_{Q_{r}}| &\leq&
C\left\{(r|\log r| + 1)\|u\|_{BMO_{p}(\Omega_{T})} +
\|u\|_{L^{1}(\Omega_{T})}\right\}\nonumber\\
&\leq& C\left(\|u\|_{BMO_{p}(\Omega_{T})} + \|u\|_{L^{1}(\Omega_{T})}\right),
\end{eqnarray}
where in the second line, we have used that $r\in (0,1)$, and that
$r|\log r|$ is bounded.\\

\noindent {\bf Step 3} (General conclusion).\\

\noindent Putting together (\ref{misang}), (\ref{result2}) and
(\ref{mwbasyd}), we get (\ref{5s5_2}). $\hfill{\blacksquare}$\\ 

\noindent We are now ready to show the proof of Proposition
\ref{prop1}.\\

\noindent {\bf Proof of Proposition \ref{prop1}.} Applying estimate
(\ref{cara_eq1}), with $m=n=1$, to the function
$\Psi\ti{u}\in W^{2,1}_{2}(\R^{2})\subseteq L^{\infty}(\R^{2})$, we get:
$$\|u\|_{L^{\infty}(\Omega_{T})}=\|\Psi\ti{u}\|_{L^{\infty}(\Omega_{T})}\leq
\|\Psi\ti{u}\|_{L^{\infty}(\R^{2})}\leq C \left(1+ \|\Psi\ti{u}\|_{BMO_{p}(\R^{2})} \left(1+
    \log^{+} \|\Psi\ti{u}\|_{W_{2}^{2,1}(\R^{2})}\right) \right).$$
Here, we have also used the fact that $\Psi=1$ over $\Omega_T$ (see
(\ref{cut_off})). Using (\ref{5s5_1}), (\ref{5s5_2}) and the above
inequality,  we directly get (\ref{prop1_eq}). $\hfill{\blacksquare}$
\subsubsection{Ideas of the proof of Theorem \ref{theorem2}}
One of the main motivations for starting with the detailed proof of
Proposition \ref{prop1} (a simplified version of Theorem \ref{theorem2})
is that it was used to show \cite[Theorem 1.1]{IJM_PI}. The other
motivation is that the arguments of the proof of Theorem \ref{theorem2}
are all contained in the proof of Proposition \ref{prop1}. It suffices
to make the following generalizations that we list below.\\

\noindent \textit{\underline{Extension of $\tilde{u}$}}. In order to extend the
function $u \in W_{2}^{2m,m}(\Omega_{T})$ to
the function $\ti{u}\in
W_{2}^{2m,m}(\widetilde{\Omega}_{T})$ with
$\widetilde{\Omega}_{T}=(-1,2)^{n}\times (-T,2T)$, we first make the
extension separately and successively with respect to the spatial
variables $x_{i}$, with $i=1\cdots n$. Then we make the
extension with respect to the time variable that is treated somehow
differently. Fix $(x_{2},.\,.\,.\,,\,x_{n}, t)\in (0,1)^{n-1}\times
(0,T)$, the spatial extension of $u$ in $x_{1}$ is as follows:
\begin{equation}\label{gene_ext}
\tilde{u}(x_{1},.\,.\,.)=
\left\{
\begin{aligned}
& \sum^{2m-1}_{j=0} c_{j}u(-\lambda_{j}x_{1},.\,.\,.)\quad &\mbox{for}&\quad
-1<x_{1}<0,\\
& \sum^{2m-1}_{j=0} c_{j}u(1+\lambda_{j}(1-x_{1}),.\,.\,.)\quad &\mbox{for}&\quad
1<x_{1}<2,
\end{aligned}
\right.
\end{equation}
with $\l_{j}=\frac{1}{2^{j}}$, and where we require that:
$$
\sum^{2m-1}_{j=0}c_{j}(-\l_{j})^{k}=1\quad \mbox{for} \quad k=0\cdots 2m-1.
$$
The above inequalities can be regarded as a linear system whose
associated matrix is of the Vandermonde type and hence invertible. This
ensures the existence of the constants $c_{j}$, $j=0\cdots 2m-1$, and
therefore the above extension (\ref{gene_ext}) gives sense.

After doing the extension with respect to $x_{1}$, the extension with
respect to $x_{2}$ is done in the same way by varying the $x_{2}$ and
fixing all other variables. This is repeated successively until the
$x_{n}$ variable.

For the time variable, we also use the same extension
(\ref{gene_ext}). Indeed, in this case, we may only sum up to $m-1$ in
(\ref{gene_ext}).\\

\noindent \textit{\underline{The cut-off function $\Psi$}}. For the definition of
the cut-off function $\Psi$, we first define the two sets:
$$\mathcal{Z}_{1}=\{(x_{1},.\,.\,.\,,\,x_{n}, t);\, \forall i=1\cdots n,\,
-1/4<x_{i}<5/4\quad
\mbox{and}\quad -T/4<t<5T/4\}$$
and
$$\mathcal{Z}_{2}=\{(x_{1},.\,.\,.\,,\,x_{n}, t);\, \forall i=1\cdots n,\,
-3/4<x_{i}<7/4\quad
\mbox{and}\quad -3T/4<t<7T/4\}.$$
The function $\Psi$ is then defined as $\Psi\in
C^{\infty}_{0}(\R^{n+1})$ with $0\leq \Psi\leq 1$ and
\begin{equation}\label{gene_cut_off}
\Psi(x,t)=
\left\{
\begin{aligned}
& 1\quad &\mbox{for}&\quad (x,t)\in \mathcal{Z}_{1}\\
& 0\quad &\mbox{for}& \quad (x,t)\in \R^{2}\setminus \mathcal{Z}_{2}.
\end{aligned}
\right.
\end{equation}

\noindent \textit{\underline{Generalization of Lemma \ref{Hgibi3}}}. An analogue
estimate of (\ref{mabaa1}) could be obtained for $(n+1)$-dimensional
parabolic cubes $Q_{r}\subseteq \widetilde{\Omega}_{T}=(-1,2)^{n}\times
(-T,2T)$. It suffices to replace $r_{2}$ satisfying (\ref{mabaa0}),
by the radius 
$$r_{n+1}=\frac{r_{n}}{\sqrt{2}},$$
where $r_{n}$ is defined recursively as follows:
$r_{j+1}=r_{j}/2$ for $0\leq j\leq n-1$.\\ 

\noindent Using the above generalizations, the proof of Theorem \ref{theorem2}
follows, line by line, the proof of Proposition
\ref{prop1}. $\hfill{\blacksquare}$
\subsection{Application of the parabolic Kozono-Taniuchi inequality}
In this subsection, we show how to apply the parabolic Kozono-Taniuchi
inequality in order to give some \textit{a priori} estimates for the
solution of certain parabolic equations. These \textit{a priori}
estimates provide a good control on the solution in order to avoid
singularities at a finite time, and hence serve for the long-time
existence. The application that will be given here deals with a model
that can be considered as a toy model. Indeed, this is a simplification of the
one treated in \cite{IJM_PI}, where a rigorous proof of the long-time
existence of solutions of a singular parabolic coupled system was presented
(see \cite[Theorem 1.1]{IJM_PI}). Consider, for $0<a<1$, the following 
parabolic equation:
\begin{equation}\label{kimi1}
\left\{
\begin{aligned}
&u_{t}(x,t)-u_{xx}(x,t)=\sin(u_{x}(x,t)u_{x}(x+a,t))+ \sin(\log
u_{x}(x,t))\quad \mbox{on}\quad \R\times(0,\infty),\\
&u(x+1,t)=u(x,t)+1\quad \mbox{on}\quad \R\times(0,\infty),\\
&u_{x}(x,0)\geq \delta_{0}>0\quad \mbox{on}\quad \R,
\end{aligned}
\right.
\end{equation}
the following proposition can be established:
\begin{proposition}{\bf (Gradient estimate)}\label{app_prop1}\\
Let $v=u_{x}$ and $m(t)=\min_{x\in\R} v(x,t)$. If $u\in
C^{\infty}(\R\times [0,\infty))$ is a smooth
solution of (\ref{kimi1}), then, for some constant
$C=C(t)>0$ we have:
\begin{equation}\label{app_eq6}
m_{t}\geq -C m \left(|\log m|+ 1\right),\quad \forall t\geq 0.
\end{equation}
\end{proposition}
\begin{rem}
Inequality (\ref{app_eq6}) directly implies that for every $t\geq 0$
we have $m(t)>0$. This is important to avoid the logarithmic
singularity in (\ref{kimi1}) when $v=u_{x}=0$.  
\end{rem}
\begin{rem}
The proof of the above proposition goes along the same lines as the
proof of \cite[Theorem 1.1]{IJM_PI}. For this reason we only present a heuristic
proof explaining only the basic ideas. The interested reader could
see the full details in \cite{IJM_PI}.
\end{rem}
{\bf Ideas of the proof of Proposition \ref{app_prop1}.} Heuristically,
the proof is divided into the following four steps. In what follows all the
constants can depend on the time $t$, but are bounded for any finite
$t$. \\

\noindent {\bf Step 1.} (First estimate from below on the gradient).\\

\noindent  Writing down the equation satisfied by $v$:
\begin{equation}\label{kimi2}
\left\{
\begin{aligned}
&v_{t}(x,t)-v_{xx}(x,t)= \cos(v(x,t)v(x+a,t))\left\{v_{x}(x,t)v(x+a,t) +
    v(x,t)v_{x}(x+a,t) \right\}\\
&\hspace{3.05cm} + \cos(\log v(x,t))
\frac{v_{x}(x,t)}{v(x,t)}\quad \mbox{on}\quad \R\times (0,\infty),\\
&v(x+1,t)=v(x,t)\quad \mbox{on}\quad \R\times (0,\infty)\\
& v(x,0)\geq \delta_{0}>0\quad \mbox{on}\quad \R,
\end{aligned}
\right.
\end{equation}
we can show that for every $t\geq 0$:
\begin{equation}\label{kimi3}
m_{t}\geq -m G\quad \mbox{with}\quad G(t)=\max_{x\in\R}|v_{x}(x,t)|.
\end{equation}

\noindent {\bf Step 2.} (Estimate of $\|v_{x}\|_{BMO_{p}}$).\\

\noindent Using the fact that $u(x+1,t)=u(x,t)+1$, and that the
right-hand term of the first equation of (\ref{kimi1}) is bounded, we apply
the $BMO$ theory for parabolic
equation to (\ref{kimi1}) and hence we obtain, for some positive
constant $c_{1}>0$:
$$
\|v_{x}\|_{BMO_{p}((0,1)\times (0,t))}\leq c_{1} \quad
\mbox{for any}\quad t>0.
$$
However, the $L^{p}$ theory for parabolic equation applied to
(\ref{kimi1}) gives, for some positive
constant $c_{2}>0$:
$$
\|v_{x}\|_{L^{1}((0,1)\times (0,t))}\leq c_{2} \quad
\mbox{for any}\quad t>0.
$$
Finally, the above two inequalities give:
\begin{equation}\label{1612}
\|v_{x}\|_{\overline{BMO}_{p}((0,1)\times (0,t))}\leq c_{1}+c_{2} \quad
\mbox{for any}\quad t>0.
\end{equation}

\noindent {\bf Step 3.} (Estimate of $\|v_{x}\|_{W^{2,1}_{2}}$).\\

\noindent Let $w=v_{x}$, we write down the equation
satisfied by $w$:
\begin{equation}\label{kimi3}
\left\{
\begin{aligned}
& w_{t}(x,t)-w_{xx}(x,t)\\
&= -\sin(v(x,t)v(x+a,t))\left\{v(x+a,t)v_{x}(x,t)+v(x,t)v_{x}(x+a,t)\right\}^{2}\\
& + \cos(v(x,t)v(x+a,t))
\left\{v(x+a,t)v_{xx}(x,t)+2v_{x}(x,t)v_{x}(x+a,t) + v(x,t)v_{xx}(x+a,t)
\right\}\\
&-\sin(\log v(x,t))\frac{v^{2}_{x}(x,t)}{v^{2}(x,t)} + \cos(\log
v(x,t))\left\{\frac{v_{xx}(x,t)}{v(x,t)} -
  \frac{v^{2}_{x}(x,t)}{v^{2}(x,t)} \right\}\quad \mbox{on}\quad
\R\times (0,\infty),\\
&w(x+1,t)=w(x,t)\quad \mbox{on}\quad \R\times (0,\infty)\\
& w(x,0)=v_{x}(x,0)\quad \mbox{on}\quad \R.
\end{aligned}
\right.
\end{equation}
Using the $L^p$ theory for parabolic equations (with various
values of $p$) to (\ref{kimi1}), (\ref{kimi2}) and (\ref{kimi3}), we
deduce, for some other positive constant $c>0$, that:
\begin{equation}\label{theeq}
\|v_{x}\|_{W^{2,1}_{2}((0,1)\times (0,t))}\leq \frac{c}{m^{2}(t)}\quad
\mbox{for any}\quad t>0.
\end{equation}
\noindent {\bf Step 4.} (Conclusion).\\

\noindent Applying the parabolic Kozono-Taniuchi inequality (\ref{prop1_eq}) to the
function $v_{x}$, using in particular (\ref{1612}) and (\ref{theeq}), we
deduce that:
$$G \leq C(1+|\log m|),$$
which, together with (\ref{kimi3}), directly gives the
result. $\hfill{\blacksquare}$\\

\noindent \textbf{Acknowledgements.} This work was supported by the
contract ANR MICA (2006-2009). The authors would like to thank M. Jazar
for his encouragement during the preparation of this work. The first
author would like to thank B. Kojok for some discussions.

\bibliographystyle{siam}
\bibliography{biblio}

\def\cprime{$'$}
\begin{thebibliography}{10}

\bibitem{BKaM}
{\sc J.~T. Beale, T.~Kato, and A.~Majda}, {\em Remarks on the breakdown of
  smooth solutions for the {$3$}-{D} {E}uler equations}, Comm. Math. Phys., 94
  (1984), pp.~61--66.

\bibitem{Bownik07}
{\sc M.~Bownik}, {\em Anisotropic {T}riebel-{L}izorkin spaces with doubling
  measures}, J. Geom. Anal., 17 (2007), pp.~387--424.

\bibitem{BreGal}
{\sc H.~Br{\'e}zis and T.~Gallou{\"e}t}, {\em Nonlinear {S}chr\"odinger
  evolution equations}, Nonlinear Anal., 4 (1980), pp.~677--681.

\bibitem{BreWai}
{\sc H.~Br{\'e}zis and S.~Wainger}, {\em A note on limiting cases of {S}obolev
  embeddings and convolution inequalities}, Comm. Partial Differential
  Equations, 5 (1980), pp.~773--789.

\bibitem{EVANS}
{\sc L.~C. Evans}, {\em Partial differential equations}, vol.~19 of Graduate
  Studies in Mathematics, American Mathematical Society, Providence, RI, 1998.

\bibitem{FraJaw90}
{\sc M.~Frazier and B.~Jawerth}, {\em A discrete transform and decompositions
  of distribution spaces}, J. Funct. Anal., 93 (1990), pp.~34--170.

\bibitem{HayWol}
{\sc N.~Hayashi and W.~von Wahl}, {\em On the global strong solutions of
  coupled {K}lein-{G}ordon-{S}chr\"odinger equations}, J. Math. Soc. Japan, 39
  (1987), pp.~489--497.

\bibitem{IJM_PI}
{\sc H.~Ibrahim, M.~Jazar, and R.~Monneau}, {\em Dynamics of dislocation
  densities in a bounded channel. {P}art {I}: smooth solutions to a singular
  coupled parabolic system}, preprint, hal-00281487.

\bibitem{CRAS_IJM}
\leavevmode\vrule height 2pt depth -1.6pt width 23pt, {\em Global existence of
  solutions to a singular parabolic/hamilton-jacobi coupled system with
  {D}irichlet conditions}, C. R. Acad. Sci. Paris, Ser. I, 346 (2008),
  pp.~945--950.

\bibitem{JohSic}
{\sc J.~Johnsen and W.~Sickel}, {\em A direct proof of {S}obolev embeddings for
  quasi-homogeneous {L}izorkin-{T}riebel spaces with mixed norms}, J. Funct.
  Spaces Appl., 5 (2007), pp.~183--198.

\bibitem{KOT03}
{\sc H.~Kozono, T.~Ogawa, and Y.~Taniuchi}, {\em Navier-{S}tokes equations in
  the {B}esov space near {$L\sp \infty$} and {BMO}}, Kyushu J. Math., 57
  (2003), pp.~303--324.

\bibitem{KT}
{\sc H.~Kozono and Y.~Taniuchi}, {\em Limiting case of the {S}obolev inequality
  in {BMO}, with application to the {E}uler equations}, Comm. Math. Phys., 214
  (2000), pp.~191--200.

\bibitem{LSU}
{\sc O.~A. Lady{\v{z}}enskaja, V.~A. Solonnikov, and N.~N. Ural{\cprime}ceva},
  {\em Linear and quasilinear equations of parabolic type}, Translated from the
  Russian by S. Smith. Translations of Mathematical Monographs, Vol. 23,
  American Mathematical Society, Providence, R.I., 1967.

\bibitem{Ogawa}
{\sc T.~Ogawa}, {\em Sharp {S}obolev inequality of logarithmic type and the
  limiting regularity condition to the harmonic heat flow}, SIAM J. Math.
  Anal., 34 (2003), pp.~1318--1330 (electronic).

\bibitem{Stockert82}
{\sc B.~St{\"o}ckert}, {\em Remarks on the interpolation of anisotropic spaces
  of {B}esov-{H}ardy-{S}obolev type}, Czechoslovak Math. J., 32 (107) (1982),
  pp.~233--244.

\bibitem{Tri92}
{\sc H.~Triebel}, {\em Theory of function spaces. {II}}, vol.~84 of Monographs
  in Mathematics, Birkh\"auser Verlag, Basel, 1992.

\end{thebibliography}
\end{document}